\documentclass[11pt]{article}
\usepackage{fullpage}

\usepackage{graphicx}
\usepackage{booktabs}
\usepackage{hyperref}

\usepackage{amsmath,amsfonts,bm}

\def\1{\bm{1}}

\def\vb{{\bm{b}}}

\def\vs{{\bm{s}}}

\def\vu{{\bm{u}}}

\def\vx{{\bm{x}}}
\def\vy{{\bm{y}}}

\def\mA{{\bm{A}}}
\def\mB{{\bm{B}}}

\def\mE{{\bm{E}}}
\def\mF{{\bm{F}}}

\def\mH{{\bm{H}}}
\def\mI{{\bm{I}}}
\def\mJ{{\bm{J}}}

\def\mW{{\bm{W}}}

\DeclareMathAlphabet{\mathsfit}{\encodingdefault}{\sfdefault}{m}{sl}
\SetMathAlphabet{\mathsfit}{bold}{\encodingdefault}{\sfdefault}{bx}{n}

\def\gH{{\mathcal{H}}}

\def\gS{{\mathcal{S}}}

\def\sR{{\mathbb{R}}}

\usepackage{nicefrac}       
\usepackage{color}          
\usepackage{algorithm}
\usepackage{algorithmic}
\usepackage{bbm, bm}
\usepackage{amsfonts, amsmath, amssymb, mathrsfs} 
\usepackage{float}
\usepackage{diagbox}
\usepackage{subcaption}
\usepackage[numbers]{natbib}
\usepackage[dvipsnames]{xcolor}
\usepackage{grffile}
\usepackage{makecell}
\usepackage{multirow}
\usepackage{wrapfig}

\newtheorem{theorem}{Theorem}
\newtheorem{definition}{Definition}
\newtheorem{lemma}{Lemma}
\newtheorem{assume}{Assumption}
\newtheorem{remark}{Remark}

\newenvironment{proof}{Proof:}{\hfill$\square$}

\newcommand{\norm}[1]{\left\|#1\right\|}

\begin{document}

\title{Explicit Superlinear Convergence Rates of Broyden's Method in Nonlinear Equations}

\author{Dachao Lin\thanks{Equal Contribution.}
	\thanks{Academy for Advanced Interdisciplinary Studies;
		Peking University;
		\texttt{lindachao@pku.edu.cn};}
	\and
	Haishan Ye\footnotemark[1]
	\thanks{School of Management; Xi'an Jiaotong University;
		\texttt{yehaishan@xjtu.edu.cn};}	
	\and
	Zhihua Zhang 
	\thanks{School of Mathematical Sciences;
		Peking University;
		\texttt{zhzhang@math.pku.edu.cn}.}
}
\maketitle
\begin{abstract}
In this paper, we study the explicit superlinear convergence rates of quasi-Newton methods. 
We particularly focus on the classical Broyden's method for solving nonlinear equations. 
We establish its explicit (local) superlinear convergence rate when the initial point is close enough to a solution and the initial Jacobian approximation is also close enough to the exact Jacobian related to the solution.
Our results present the explicit superlinear convergence rates of Broyden's ``good'' and ``bad'' update schemes \cite{martinez2000practical}.
These explicit convergence rates in turn provide some important insights on the performance difference between the ``good'' and ``bad'' schemes, which is also validated empirically.

\end{abstract}

\section{Introduction}

We are concerned with the following nonlinear equation system
\begin{equation}\label{obj}
    \mF(\vx) = \bm{0},
\end{equation}
where $\mF: \sR^n \to \sR^n$ is a differentiable vector-valued function with $\mF(\vx) = (F_1(\vx), \dots, F_n(\vx))^\top, \vx \in \sR^n$. It is an important tool in describing key properties of control systems \cite{more1989collection}, economic systems \cite{albu2006non}, stationary points of optimization problems, etc.
Usually, it is hard to obtain an analytic solution to  nonlinear equations.
Thus, many numerical methods have been proposed to solve  nonlinear equations.

Newton's method is a classical approach for this problem based on the following iteration
\begin{align*}
	\vx_{k+1} = \vx_{k} - \left[\mJ(\vx_{k})\right]^{-1}\mF(\vx_{k}),
\end{align*}
where $\mJ(\vx)\in \sR^{n\times n}$ is the exact Jabobian of $\mF$ at point $\vx$. 
The above iteration converges under some conditions, essentially requiring $\vx_0$ to be sufficiently close to a solution.
The convergence is fast, usually quadratic, under some mild conditions \cite[Chapter 11.1]{nocedal2006numerical}. 
Due to the fast convergence rate, Newton's method plays an important tool in solving nonlinear equations and lies at the heart of many important algorithms.
However,  Newton's method suffers from expensive computation  because it requires computing the exact Jacobian which is $n^2$ first order derivatives, and the inverse of the Jacobian which takes $O(n^3)$ flops.

To remedy the heavy computational costs of Newton's method, quasi-Newton methods have been proposed, which do not need to compute the Jacobian and its inverse in each iteration.
Rather than computing the exact Jacobian, the quasi-Newton method attempts to give an approximation to the Jacobian and update it at each iteration so that it mimics the behavior of the true Jacobian $\mJ(\vx)$.
For each iteration of the quasi-Newton method, the main computational cost is to compute $\mF(\vx)$.
Therefore, the quasi-Newton method can also achieve a fast convergence rate similar to  Newton's method but at a low computational cost for each iteration. 
There have been several quasi-Newton versions proposed, such as Broyden's method \cite{broyden1965class}, ABS \cite{abaffy1989abs}, etc.
Among quasi-Newton methods in solving nonlinear equations, Broyden's method is the most important one.
Thus, in this paper, we focus on Broyden's method for solving nonlinear equations, including Broyden's ``good'' scheme and Broyden's ``bad'' scheme \cite{al2014broyden,dennis1996numerical,griewank2012broyden,martinez2000practical}.

Because of the advantages in solving nonlinear equations, Broyden's method is a hot research topic and many improvements over Broyden's method have been made to achieve better performance \cite{sun1997newton,lukvsan2017new,schubert1970modification}.
At the same time, the convergence analysis of Broyden's method attracts much attention to figure out its exact convergence properties.
Many works have shown that Broyden's method can achieve a superlinear convergence rate which theoretically validates its fast convergence rate \cite{broyden1973local,kelley1991new,dennis1996numerical}.
Whereas, the superlinear convergence rates provided in previous works are only asymptotic, that is,
these results only show that the ratio of successive residuals tends to zero, i.e., 
$\lim_{k\to +\infty} \norm{\bm{x}_{k+1}-\bm{x}_{*}} / \norm{\bm{x}_k-\bm{x}_{*}} = 0$,
where $\{\vx_k\}_{k=1}^{+\infty}$ is the iterative update sequence, and $\vx_*$ is a nondegenerate solution (i.e., the Jacobian $\mJ(\vx_*)$ is nonsingular).
However, an asymptotic superlinear convergence rate can not describe the exact convergence speed.  
For example, it is unknown whether the residuals converge
like $O(c^{k^2})$, or $O(k^{-k})$, where $c \in (0,1) $ is some constant \cite{rodomanov2021greedy}.
And the convergence speeds described by these superlinear rates vary greatly.
Thus, an explicit superlinear rate for Broyden's method is desirable for obtaining the convergence dynamics of algorithms and  judging the performance of algorithms in theory. 
Moreover, without the explicit superlinear rate, it is hard to figure out the factors which mainly determine the performance of algorithms.
For example, the reason for the performance difference between Broyden's ``good'' scheme and ``bad'' scheme is an important open problem in the study of quasi-Newton methods for solving nonlinear equations \cite{martinez2000practical}.


In this paper, we aim to give the explicit superlinear convergence rates of Broyden's method including the ``good'' and ``bad'' schemes, and provide some insights into the reason for the performance difference between the ``good'' and ``bad''  in terms of the explicit rates. 
We summarize our contribution as follows:
\begin{enumerate}
    \item We first show that if the initial point $\vx_0$ and the initial Jacobian approximation $\mB_0$ are sufficiently close to one solution $\vx_*$ and corresponding $\mJ_*:=\mJ(\vx_*)$, then the sequence produced by Broyden's ``good'' scheme converge to $\vx_*$ at a superlinear rate of the form (shown in Theorem \ref{thm:good}):
    \[ \left( \frac{\norm{\mJ_*^{-1}\mB_0 - \mI_n}_F + \sqrt{\|\mJ_*^{-1}\| \cdot \norm{\vx_0 - \vx_*}}}{\sqrt{k}} \right)^{k}. \]
	Moreover, we also show that Broyden's ``bad'' scheme converges with a superlienar rate of the form (shown in Theorem \ref{thm:bad}):
	\[ \left( \frac{\norm{\mJ_*\mB_0^{-1} - \mI_n}_F + \sqrt{\|\mJ_*^{-1}\|^2 \cdot \norm{\mJ_*\left(\vx_0 - \vx_*\right)}}}{\sqrt{k}} \right)^{k}. \]
	\item Our explicit superlinear rates of the ``good'' and ``bad'' schemes show that two initial distances, i.e., $\norm{\mJ_*^{-1}\mB_0 - \mI_n}_F$ (or $\norm{\mJ_*\mB_0^{-1} - \mI_n}_F$), the distance between the initial Jacobian and the exact Jacobian, and $\norm{\vx_0-\vx_*}$ (or $\norm{\mJ_*(\vx_0-\vx_*)}$), the distance between the initial point and the solution, determine the convergence speed of Broyden's method. 
	Meanwhile, we also gives the exact description of the neighborhood around the solution to guarantee the convergence of Broyden's method and show a natural trade-off between the initial conditions on $\vx_0, \mB_0$ and the superlinear convergence rate. 
	Based on our convergent results, we find that Broyden's ``bad'' scheme prefers to an overestimated $\mB_0 =\alpha \mJ_*, \alpha \geq 1$, while the ``good'' scheme is suitable to a underestimated $\mB_0 =\alpha \mJ_*, 0 < \alpha < 1$.
\end{enumerate}

\subsection{Organization}
The remainder of the paper is organized as follows. After having reviewed related work in Section \ref{sec:related}, we present some
preliminaries in Section \ref{sec:preliminaries}, including the notation and introduction of Broyden's ``good'' and ``bad'' update schemes for obtaining an approximate Jacobian matrix.
In Section \ref{sec:good-bro}, we analyze both Broyden's ``good'' and ``bad'' schemes and present our main results of explicit local superlinear convergence rates.
In Section \ref{sec:exp}, we perform numerical experiments to support our theoretical results.
We also show some comparison with previous work in Section \ref{sec:compare}.
Finally, we conclude our work in Section \ref{sec:conclusion}.

\section{Related Work}\label{sec:related}

The quasi-Newton methods for solving nonlinear equations have been widely studied and several important algorithms have been proposed \cite{broyden1965class,thomas1975sequential,abaffy1989abs,spedicato1978some}.
Among all the quasi-Newton algorithms for solving nonlinear equations, Broyden's method is the most famous and important \cite{broyden1965class,broyden1967quasi,kelley1991new,nocedal2006numerical}.
For each iteration, Broyden's method conducts a rank-one update for the approximate Jacobian.
The updated matrix has to satisfy so-called secant equations.
According to the way of updating, there exist  Broyden's ``good'' scheme, Broyden's ``bad'' scheme and Broyden's symmetric update (identical to the famous SR1 method) \cite{broyden1965class}.
Broyden's ``good'' scheme is probably still the most common choice.
In contrast, the ``bad'' scheme may perform not so well but on a small set of test problems \cite{al2014broyden}. The reasons why Broyden's ``good'' scheme is good and Broyden's ``bad'' scheme is bad are not well understood \cite{martinez2000practical}.

The quasi-Newton variants of \cite{spedicato1978some} were derived from the optimal conditioning criterion proposed in \cite{oren1976optimal}.
For the unconstrained optimization, the quasi-Newton is also an important class of algorithms which includes the famous Broyden-Fletcher-Goldfarb-Shanno (BFGS) method \cite{broyden1970convergence2,broyden1970convergence,fletcher1970new,goldfarb1970family,shanno1970conditioning}, the Davidon-Fletcher-Powell (DFP) method \cite{davidon1991variable,fletcher1963rapidly} and the Symmetric Rank 1 (SR1) method \cite{broyden1965class,broyden1967quasi,davidon1991variable}. 
Though there are some important differences between quasi-Newton methods for the unconstrained optimization and for solving nonlinear equations, the most widely used quasi-Newton methods can typically achieve superlinear convergence rates.  
\cite{powell1971convergence} showed that DFP can obtain the local superlinear convergence.
Since then, many works have provided asymptotic superlinear convergence results of quasi-Newton methods for the unconstrained optimization \cite{byrd1987global,griewank1982local,kovalev2020fast,stachurski1981superlinear,yabe1996local}.
The work \cite{broyden1973local} first proved that Broyden's method for solving nonlinear equations converges superlinearly.
After that, many works gave improved analysis for Broyden's method \cite{dennis1996numerical,gruver1981algorithmic,kelley1991new,hwang1992convergence,kelley1995iterative}. 

Recently, \cite{rodomanov2021greedy} gave the first explicit superlinear convergence rate for the greedy quasi-Newton method which takes updates by greedily selecting from basis vectors to maximize a certain measure of progress. 
This work establishes an explicit non-asymptotic rate of the local superlinear convergence $(1-\frac{1}{n\varkappa})^{k^2/2}$, where $k$ is the iteration number and $\varkappa$ is the condition number of the objective function in question. 
\cite{lin2021faster} improved the rates of greedy quasi-Newton updates as well as random quasi-Newton updates,  obtaining a faster condition-number-free superlinear convergence rate $(1-\frac{1}{n})^{k^2/2}$. 
As for classical quasi-Newton methods, \cite{rodomanov2021rates} also analyzed the classical well-known DFP and BFGS methods, adopting standard Hessian update direction through the previous variation. 
They demonstrated the rates of the forms $(\frac{n\varkappa^2}{k})^{k/2}$ and $(\frac{n\varkappa}{k})^{k/2}$ for the standard DFP and BFGS methods, respectively. 
Such rates have faster initial convergence rates, while slower final rates compared to \cite{rodomanov2021greedy}'s results.
Furthermore, \cite{rodomanov2021new} also improved the results of \cite{rodomanov2021rates} to $(\frac{n\varkappa \ln \varkappa}{k})^{k/2}$ and $(\frac{n\ln\varkappa}{k})^{k/2}$ for the standard DFP and BFGS methods, though having similar worse long-history behavior compared to \cite{rodomanov2021greedy,lin2021faster}.
\cite{ye2021explicit} extended the results of \cite{rodomanov2021new} to the modified SR1 method with the correction strategy.

The above results \cite{rodomanov2021greedy,rodomanov2021new,rodomanov2021rates,lin2021faster,ye2021explicit} only make the assumption that the initial points should be in a small region near the unique solution, that is, $\norm{\vx_0-\vx_*}$ should be small enough.
Based on another assumption that the initial approximate Hessian is close enough to the exact Hessian, \cite{jin2020non} provided a dimension-free superlinear convergence rate for the standard BFGS and DFP methods. 
The assumption used in \cite{jin2020non} is much similar to ours. 
Our work in this paper requires the assumption that both $\norm{\vx_0-\vx_*}$ and $\norm{\mB_0-\mJ(\vx_*)}$ should be small enough.
Such requirements is unavoidable due to the non-singularity of approximate Jacobian matrix based on Broyden’s update. 
Meanwhile, our work is not a simple extension of \cite{jin2020non} because Broyden's method for solving nonlinear equations is much different from DFP and BFGS (see Section \ref{sec:compare} for detail). 

To our best knowledge, the result about the explicit local superlinear convergence rate of the original Broyden's method for solving general nonlinear functions is still unknown. \cite{griewank1987local} obtained a similar superlinear convergence rate $(\frac{O(1)}{k})^{k/2}$  as ours but for a modified Broyden's scheme. 
Their proof is mainly applied to a line search version of the classical Broyden's ``good'' scheme, which actually makes the analysis much simpler. 
Moreover, they did not explicitly show the relationship between the initial conditions and the superlinear convergence rates because they hid such factors to unknown constants. 
Furthermore, \cite{griewank1987local} lacked the explicit superlinear rate of Broyden's ``bad'' scheme.
Therefore, our result in this paper is novel to the analysis of Broyden's method.

\section{Preliminaries}\label{sec:preliminaries}
\subsection{Notation}
We denote vectors by lowercase bold letters (e.g., $ \vu, \vx$), and matrices by capital bold letters (e.g., $ \mW = [w_{i j}] $).
We use $[n]:=\{1,\dots,n\}$ and $\mI_n$ is the $\sR^{n \times n}$ identity matrix, and $\mathrm{Unif}(\gS^{d-1})$ as the uniform distribution from $\gS^{d-1}$.
Moreover, $\norm{\cdot}$ denotes the $\ell_2$-norm (standard Euclidean norm) for vectors, or spectral norm for a given matrix: $\norm{\mA} = \sup_{\norm{\vu}=1, \vu\in\sR^n}\norm{\mA\vu}$, and $\norm{\cdot}_F$ denotes the Frobenius norm of a given matrix: $\norm{\mA}_F = \sqrt{\sum_{i=1}^{m} \sum_{j=1}^n a_{ij}^2}$, where $\mA =[a_{i j}] \in \sR^{m\times n}$.
For two real $n\times n$ symmetric matrices $\mA$ and $\mB \in \sR^{n \times n}$, we denote $\mA\succeq\mB$ if $\mA-\mB$ is a positive semidefinite matrix.
We use the standard $O(\cdot), \Omega(\cdot)$ and $\Theta(\cdot)$ notation to hide universal constant factors.

\subsection{Broyden's Update}

\begin{algorithm}[t]
	\caption{Broyden's ``good'' scheme (only for convergence analysis)}
	\begin{algorithmic}[1]
		\STATE Initialization: set $\mB_0, \vx_0$.
		\FOR{$ k \geq 0 $}
		\STATE Update $\vx_{k+1} = \vx_{k} - \mB_k^{-1} \mF(\vx_k)$.
		\STATE Set $\vu_{k}=\vx_{k+1}-\vx_k$, $\vy_{k}= \mF(\vx_{k+1})- \mF(\vx_{k})$.
		\STATE Compute 
		$\mB_{k+1} =  \mB_k + \frac{\left(\vy_k-\mB_k\vu_k\right)\vu_k^\top}{\vu_k^\top\vu_k}$.
		\ENDFOR
	\end{algorithmic}
	\label{algo:broyden}
\end{algorithm}

Before starting our theoretical results, we briefly review Broyden's update. We focus on solving a nonlinear equation system defined in Eq.~\eqref{obj}. 
Broyden's method iteratively updates the approximate Jacobian $\{\mB_k\}_{k=0}^{+\infty}$ to approach the true Jacobian matrix, while updating the parameters $\{\vx_k\}_{k=0}^{+\infty}$ as:
\begin{equation}\label{eq:x-update}
    \vx_{k+1} = \vx_k-\mB_k^{-1}\mF(\vx_k).
\end{equation}
At the same time, the alternative matrix $\mB_k$ needs to satisfy the secant condition, that is 
\begin{equation}\label{eq:secant}
\mB_k\vu_k=\vy_k, \mbox{ with }\vu_k = \vx_{k+1}-\vx_k,\; \vy_{k}=\mF(\vx_{k+1})-\mF(\vx_k).
\end{equation} 
However, we could not obtain a well-definite solution with only Eq.~\eqref{eq:secant} because it is an underdetermined linear system with an infinite number of solutions. 
Previous works adopt several different extra conditions to give a reasonable solution. 
The most common concern chooses the Jacobian approximation sequence $\mB_{k+1}$, which is as close as possible to $\mB_{k}$ under a measure. Broyden's ``good'' update \cite{broyden1965class,kelley1995iterative} is defined as follows 
\begin{equation}\label{eq:B-update-good}
	\mB_{k+1} =  \mB_k + \frac{\left(\vy_k-\mB_k\vu_k\right)\vu_k^\top}{\vu_k^\top\vu_k},
\end{equation}
which is obtained by solving the following problem:
\begin{equation*}
\min_{\mB} \norm{\mB-\mB_k}_F, \ s.t., \mB\vu_k=\vy_k.
\end{equation*} 
The detailed derivation of Broyden's ``good'' scheme can be found in  Theorem 4.1 in \cite{dennis1977quasi} and the detailed algorithmic description is listed in Algorithm~\ref{algo:broyden}.
Note that, Algorithm~\ref{algo:broyden} is \textit{only for convergence analysis}. 
For the practical implementation, researchers would employ the Sherman-Morrison-Woodbury formula \cite{sherman1950adjustment} for updating the inverse of matrix $\mH_k := \mB_{k}^{-1}$ directly as follows:
\[ \mH_{k+1} = \mH_k-\frac{\left(\mH_k\vy_k - \vu_k\right) \vu_k^\top\mH_k}{\vu_k^\top\mH_k\vy_k}. \]
Similarly, Broden's ``bad'' scheme \cite{broyden1965class,martinez2000practical} is defined as follows
\begin{equation*}
    \mH_{k+1} = \mH_k + \frac{\left(\vu_k-\mH_k\vy_k\right)\vy_k^\top}{\vy_k^\top\vy_k},
\end{equation*}
which is obtained by solving the following problem:
\[ \min_{\mH} \norm{\mH-\mH_k}_F, \ s.t., \mH\vy_k = \vu_k. \]
The update of Broyden's ``bad'' scheme tries to directly approximate the inverse of the exact Jacobian under the constraint of the secant equation.
Moreover, we could obtain its $\mB_k$ update scheme as follows:
\[ \mB_{k+1} = \mB_k-\frac{\left(\mB_k\vu_k - \vy_k\right) \vy_k^\top\mB_k}{\vy_k^\top\mB_k\vu_k}. \]
The detailed implementation of Broyden's ``bad'' scheme is described in Algorithm~\ref{algo:broyden-bad}. 
Moreover, we need to clarify that Broyden \cite{broyden1965class} suggested this formula as well, but it seems that he and others had less favourable numerical experience, which lead to the moniker ``bad'' Broyden update \cite{griewank2012broyden,dennis1996numerical}.
However, it is not clear that, in practice, the ``good'' scheme is really better than the ``bad'' scheme \cite{martinez2000practical}. Indeed, sometimes, the ``bad'' scheme is not at all bad \cite{kvaalen1991faster}.

For better understanding the update rules of Broyden's method, we turn to a simple linear objective $\mF(\vx) =\mA\vx-\vb$ for explanation. 
When applied to linear objective, we could simplify the update rules using the fact $\vy=\mA\vu$. 
Thus, we could simplify and unify  Broyden's ``good'' and ``bad'' updates  below.
\begin{definition}[Broyden's Update]
	Letting $\mB, \mA\in\sR^{n \times n}$ and $\vu\in\sR^n$,  we define
    \[ \text{Broyd} \left(\mB, \mA, \vu\right) := \mB +  \frac{\left(\mA-\mB\right)\vu\vu^\top}{\vu^\top\vu}. \]
\end{definition}
Then, for the simple linear objective $\mF(\vx) =\mA\vx-\vb$, Broyden's ``good'' update can be written as $\mB_{k+1}=\text{Broyd} \left(\mB_k, \mA, \vu_k\right)$, and the ``bad'' one is written as $\mH_{k+1}=\text{Broyd} \left(\mH_k, \mA^{-1}, \vy_k\right)$. 

\begin{algorithm}[t]
	\caption{Broyden's ``bad'' scheme (only for convergence analysis)}
	\begin{algorithmic}[1]
		\STATE Initialization: set $\mH_0, \vx_0$.
		\FOR{$ k \geq 0 $}
		\STATE Update $\vx_{k+1} = \vx_{k} - \mH_k \mF(\vx_k)$.
		\STATE Set $\vu_{k}=\vx_{k+1}-\vx_k$, $\vy_{k}= \mF(\vx_{k+1})- \mF(\vx_{k})$.
		\STATE Compute $\mH_{k+1} = \mH_k + \frac{\left(\vu_k - \mH_k \vy_k\right)\vy_k^\top}{\vy_k^\top\vy_k}$.
		\ENDFOR
	\end{algorithmic}
	\label{algo:broyden-bad}
\end{algorithm}

\subsection{Notation for Convergence Analysis}
To give a better understanding of update, we will need the integral Jacobian between two adjacent points following Algorithms \ref{algo:broyden} and \ref{algo:broyden-bad}:
\begin{equation}\label{eq:int-J}
    \mJ_k := \int_{0}^{1} \mJ\left(\vx_k+t\vu_k\right)d t, \ \vu_k := \vx_{k+1}-\vx_k.
\end{equation}
To estimate the convergence rate of Algorithm~\ref{algo:broyden}, we let the norm between $\vx_k$ and a solution $\vx_*$ as
\begin{equation}\label{eq:r}
	r_k := \norm{\vx_k-\vx_*}.
\end{equation}
Moreover, we introduce the potential function that measures the objective gap between $\vx_k$ and $\vx_*$ as follows:
\begin{equation}\label{eq:lam}
    f_{*}(\vx) := \norm{ \mJ_*^{-1}\left[\mF(\vx)-\mF(\vx_*)\right] } = \norm{ \mJ_*^{-1}\mF(\vx) }, \mbox{ and } f_k := f_{*}(\vx_k).
\end{equation}
The explicit convergence rate of Algorithm~\ref{algo:broyden} can be described by $f_k$.
We also introduce a potential function which measures the approximation precision of the Jacobian as follows 
\begin{equation}\label{eq:jac-measure}
    \overline{\mB}:= \mJ_*^{-1}\left(\mB-\mJ_{*}\right)= \mJ_*^{-1}\mB-\mI_n, \ \sigma_{*}(\mB) := \norm{\overline{\mB}}_F, \mbox{ and } \sigma_k := \sigma_{*}(\mB_k).
\end{equation}
When analyzing Algorithm~\ref{algo:broyden-bad}, we adopt another group of measures.
We use the scaled distance
\begin{equation}\label{eq:R}
	R_k := \norm{\mJ_*\left(\vx_k-\vx_*\right)}
\end{equation}
to replace $r_k$.
In contrast to using $\sigma_{k}$ to measure the distance of the approximate Jacobian to the exact one, we introduce $\tau_{k}$  defined as follows:
\begin{equation}\label{jac-measure-bad}
	\widetilde{\mH} := \mJ_*\left(\mH - \mJ_*^{-1}\right) = \mJ_*\mH-\mI_n, \; \tau_{*}(\mH) := \|\widetilde{\mH}\|_F, \mbox{ and } \tau_k := \tau_{*}(\mH_k),
\end{equation}
which measures the distance between the approximate Jacobian inverse and the exact one. 
Moreover, we directly use 
\begin{equation}\label{eq:fk}
    F_*(\vx):=\norm{\mF(\vx)-\mF(\vx_*)}, \mbox{ and } F_k:=F_*(\vx_k)
\end{equation}
instead of $f_k$ to derive the explicit convergence rate of Algorithm~\ref{algo:broyden-bad}.
Finally, we also denote $\mu := 1/\|\mJ_*^{-1}\|, L := \|\mJ_*\| $, and $ \varkappa := L / \mu$ as the condition number in our analysis.
Then we have
\begin{equation}\label{eq:rR}
    \mu r_k \leq R_k \leq L r_k.
\end{equation}
\subsection{The Assumptions}

Before analyzing the convergence rate of Broyden's method, we first assume that the Jacobian is $M$-Lipschitz continuous related to $\vx_*$ which is described as follows.
\begin{assume}\label{ass:lisp}
	The Jacobian $\mJ(\vx)$ is Lipschitz continuous related to $\vx_*$ with
	parameter $M$ under the spectral norm, i.e., 
	\begin{equation}\label{eq:lisp}
	    \norm{\mJ(\vx)-\mJ(\vx_*)} \leq M \norm{ \vx-\vx_* }, \forall \vx \in \sR^n.
	\end{equation}
\end{assume}
The assumption of Lipschitz continuity is a standard assumption in analyzing the convergence rate of Newton's method and Broyden's method for solving nonlinear equations \cite{nocedal2006numerical,kelley1995iterative,dennis1996numerical}.
Particularly, Assumption \ref{ass:lisp} has appeared in \cite[Lemma 8.2.1]{dennis1996numerical}, and a similar assumption for Hessian matrix also appears in \cite[Assumption 6.2]{nocedal2006numerical} and \cite[Assumption 3.2]{jin2020non}, which analyzes the convergence rate of the quasi-Newton methods for unconstrained optimization.  
Note that, Assumption~\ref{ass:lisp} is weaker than the strong self-concordance recently proposed by \cite{rodomanov2021greedy,rodomanov2021new,rodomanov2021rates}.
Indeed, our analysis only requires that Assumption \ref{ass:lisp} holds near the optimal solution $\vx_*$. For brevity, we assume Eq.~\eqref{eq:lisp} holds in $\sR^n$.

Now we introduce another assumption that is important in analyzing the convergence rate of Broyden's method.
\begin{assume}\label{ass:optimal}
	Sequences $\{\vx_{k}\}_{k=0}^{+\infty}$ and $\{\mB_{k}\}_{k=0}^{+\infty}$ generated by Broyden's methods (including the ``good'' and ``bad'' schemes) are well-defined, that is, $\mB_k$ is nonsingular for all
	$k$, and $ \vx_{k+1} $ and $ \mB_{k+1} $ can be  computed from $ \vx_{k} $ and $ \mB_{k}$. Moreover, the solution $\vx_*$ of Eq.~\eqref{obj} is nondegenerate, i.e., $\mJ_{*}:=\mJ(\vx_*)$ is nonsingular.
\end{assume}
Assumption~\ref{ass:optimal} is a standard assumption in analyzing the convergence rate of Broyden's method which also appears in \cite{dennis1996numerical,nocedal2006numerical}. 
It assumes that the sequence $ \{\mB_k\}_{k=0}^{+\infty} $ of Broyden's updates exists and keeps nonsingular forever. 
Indeed, if $\vx_0$ is close enough to a nondegenerate solution $\vx_*$ and $\mJ_0$ is also close enough to $\mJ_*$, we could derive that $\mB_k$ (or $\mH_k$) is nonsingular for all $k\geq 0$ under Assumption \ref{ass:lisp}.

\section{Main Results}
\label{sec:good-bro}

In this section, we try to provide explicit superlinear convergence rates of Broyden's method including the ``good'' and ``bad'' schemes.
It is well-known and just as mentioned in \cite[Chapter 11.1]{nocedal2006numerical} that when the starting point is remote from a solution, the behavior of Newton’s method can be erratic.
Hence, our goal is to give explicit superlinear convergence provided that the method can start from an initial point $\bm{x}_0$, which is close enough to a nondegenerate solution $\vx_*$. 
In the next, we will first provide the explicit superlinear rates of Broyden's ``good'' and ``bad'' schemes.
Then we will compare the rates from different perspectives and provide some insights on the performance difference between the ``good'' and ``bad'' schemes.

\subsection{Convergence Rate of Broyden's ``Good'' Scheme}
\label{subsec:good}

By Eq.~\eqref{eq:x-update} and Eq.~\eqref{eq:B-update-good}, we have the following updates of Broyden's ``good'' scheme: 
\begin{equation}\label{eq:x-up}
\left\{
\begin{aligned}
& \vx_{k+1} = \vx_k-\mB_k^{-1} \mF(\vx_k), \\
& \mB_{k+1} = \mB_k +  \frac{\left(\vy_k - \mB_k \vu_k\right)\vu_k^\top}{\vu_k^\top\vu_k},
\end{aligned}
\right. 
\end{equation}
where $\vy_k = \mF(\vx_{k+1})-\mF(\vx_k)$ and $\vu_k=\vx_{k+1}-\vx_k$.
We then have $\vy_k\stackrel{\eqref{eq:int-J}}{=}\mJ_k\vu_k$.
Accordingly, the update of the Jacobian approximation can be represented as
\begin{equation*}
\mB_{k+1} = \mB_k + \frac{\left(\mJ_k - \mB_k \right)\vu_k\vu_k^\top}{\vu_k^\top\vu_k} = \text{Broyd}\left(\mB_k, \mJ_k, \vu_k\right).
\end{equation*} 

First, we give a lemma, showing that $r_k$ can achieve a local linear convergence rate under benign conditions, i.e., the starting point $\vx_0$ and $\mB_0$ are close enough to $\vx_*$ and $\mJ_*$, respectively.

\begin{lemma}\label{lem:base}
	Suppose the objective function $\mF(\cdot)$ satisfies Assumption~\ref{ass:lisp} and sequences $\{\vx_k\}_{k=0}^{+\infty}$ and $\{\mB_k\}_{k=0}^{+\infty}$ generated by Algorithm~\ref{algo:broyden} satisfy Assumption~\ref{ass:optimal}. 
	If there exists a $q \in (0, 1)$ such that $r_0$ and $\sigma_0$ satisfy
	\begin{equation}\label{eq:init-cond}
		\sigma_0 \leq \frac{q}{q+1}, \quad \frac{M r_0}{\mu} \leq \frac{q(1-q)}{8} \left(\frac{q}{1+q}-\sigma_0\right),
	\end{equation}
	then for all $k\geq 0$ it holds that 
	\begin{equation}\label{eq:sigma-bound}
		\sigma_k^2 \leq \sigma_0^2 + \frac{1+q}{1-q}\left( \frac{M r_0}{\mu} + \frac{M^2r_0^2}{\mu^2} \right) \leq \left(\frac{q}{q+1}\right)^2,
	\end{equation}
	and 
	\begin{equation}\label{eq:linear-con}
		r_{k+1} \leq q r_k.
	\end{equation}
\end{lemma}

Based on the linear convergence of $r_k$, we are ready to prove the local superlinear convergence rate of Broyden's ``good'' scheme.

\begin{lemma}\label{lem:super}
	Suppose the objective function $\mF(\cdot)$ and sequences $\{\vx_k\}_{k=0}^{+\infty}$ and $\{\mB_k\}_{k=0}^{+\infty}$ satisfy the properties described in Lemma~\ref{lem:base}.  
	Let distances $r_0$ and $\sigma_0$ satisfy Eq.~\eqref{eq:init-cond} for some $q \in (0, 1)$.
	Then, for all $k\geq 1$, it holds that  
    \[ f_{k} \leq \left[\frac{q^2}{k}\right]^{k/2}f_0. \]
\end{lemma}
\begin{remark}
	Lemmas~\ref{lem:base} and~\ref{lem:super} show that when the initial points $\vx_0$ and $\mB_0$ are close enough to one solution $\vx_*$ and its  corresponding Jacobian $\mJ_*$, respectively, $r_k$ converges with a linear rate and $f_k$ converges superlinearly. 
	Moreover, these lemmas also reveal that the convergence rate $q$ depends on the radius of the neighborhood at initialization and there is a trade-off between the convergence speed and the radius, that is,  
	the convergence speed could be faster when a smaller region around a solution at initialization is provided and vice versa. 
	Thus, one can improve the superlinear convergence rates at the cost
	of reducing the radius of the neighborhood at initialization.
	Similar point has been discovered in the analysis of the quasi-Newton method for the unconstrained optimization in \cite{jin2020non}.
	However, such initialization is hard or expensive to satisfy generally.
	Furthermore, the distances of $\norm{\vx_0-\vx_*}$ and $\norm{\mB_0-\mJ_{*}}$ also have relevance given a fixed convergence rate $q$. 
\end{remark}

To better quantify the superlinear convergence rate of Broyden's ``good'' scheme, we define $q_m$ given the suitable $\sigma_0$ and $r_0$:
\begin{equation*}
	q_m := \min \; q, \ s.t. \ \sigma_0 \leq \frac{q}{q+1}, \ \frac{M r_0}{\mu} \leq \frac{q(1-q)}{8} \left(\frac{q}{1+q}-\sigma_0\right).
\end{equation*}
Such a $q_m$ needs a pair of proper $\sigma_0, r_0$ to guarantee nonempty feasible set, and $q_m$ determines the optimal superlinear convergence rate of Broyden's ``good'' scheme from Lemmas~\ref{lem:base} and~\ref{lem:super}.

\begin{theorem}\label{thm:good}
	Suppose the objective function $\mF(\cdot)$ satisfies Assumption~\ref{ass:lisp} and sequences $\{\vx_k\}_{k=0}^{+\infty}$ and $\{\mB_k\}_{k=0}^{+\infty}$ generated by Algorithm~\ref{algo:broyden} satisfy Assumption~\ref{ass:optimal}.  
	If $r_0$ and $\sigma_0$ satisfy
	\begin{equation}\label{eq:init-super}
		\frac{32 M r_0}{\mu} + \sigma_0 \leq \frac{1}{3},
	\end{equation}
    then it holds that $q_m$ is well-defined and  $q_m = \Theta\left(\sigma_0+\sqrt{\frac{M r_0}{\mu}}\right)$, 
	and we have
	\begin{equation}\label{eq:easyread-convergence}
		\begin{aligned}
			r_{k+1} & \leq 6\left(\sigma_0+\sqrt{\frac{M r_0}{\mu}}\right)r_k, \forall k\geq 0, \\
			f_{k} & \leq \left[\frac{6\left(\sigma_0+\sqrt{\frac{M r_0}{\mu}}\right)}{\sqrt{k}}\right]^{k}f_0, \forall k\geq 1.
		\end{aligned}
	\end{equation}
\end{theorem}
\begin{proof}
	To give an upper and lower bound of $q_m$, we consider 
	\[ f(q) := q(1-q)\left(\frac{q}{1+q}-\sigma_0\right), \quad q \geq \frac{\sigma_0}{1-\sigma_0}. \]
	We only need to consider the second constraint in the definition of $q_m$. 
	
	Note that $q(1-q)$ is increasing when $q \leq \frac{1}{2}$, and $\frac{q}{q+1}$ is also 
	increasing when $q>0$. 
	Hence $f(q)$ is increasing when $\frac{1}{2} \geq q \geq \frac{\sigma_0}{1-\sigma_0}$, where we use $\frac{1}{2} \geq \frac{\sigma_0}{1-\sigma_0}$ by Eq.~\eqref{eq:init-super}. 
	Moreover, we have
	\begin{equation}\label{eq:f-2}
	    f\left(\frac{1}{2}\right) = \frac{1}{4} \left(\frac{1}{3}-\sigma_0\right) \stackrel{\eqref{eq:init-super}}{\geq} \frac{8 M r_0}{\mu},
	\end{equation}
	showing that $q=\frac{1}{2}$ is in the feasible set. Thus $q_m$ is well-defined and 
	\begin{equation}\label{eq:qm-2}
	     0 \leq \sigma_0 \leq \frac{\sigma_0}{1-\sigma_0} \leq q_m \leq \frac{1}{2}.
	\end{equation} 
	
	Since $f(\frac{\sigma_0}{1-\sigma_0}) = 0 \leq \frac{8Mr_0}{\mu} \stackrel{\eqref{eq:f-2}}{\leq} f(\frac{1}{2})$, we get $f(q_m) = \frac{8Mr_0}{\mu}$.
	Then we obtain
	\[ \frac{8 M r_0}{\mu} = f(q_m) \stackrel{\eqref{eq:qm-2}}{\geq} \frac{q_m}{2}\left(\frac{q_m}{2}-\sigma_0\right), \]
	which gives
	\begin{equation}\label{eq:qm-upper}
		q_m \stackrel{\eqref{eq:qm-2}}{\leq} \sigma_0 +\sqrt{\sigma_0^2+\frac{32 M r_0}{\mu}} \leq 2\sigma_0+6\sqrt{\frac{M r_0}{\mu}} \leq 6\left(\sigma_0+\sqrt{\frac{M r_0}{\mu}}\right).
	\end{equation}
	Furthermore, from $q_m \geq 0$ we have that
    \[ \frac{8 M r_0}{\mu}= f(q_m) \stackrel{\eqref{eq:qm-2}}{\geq} q_m\left(q_m-\sigma_0\right), q_m \geq \sigma_0, \]
	which gives
	\begin{equation}\label{eq:qm-lower}
		q_m \stackrel{\eqref{eq:qm-2}}{\geq}  \frac{\sigma_0+ \sqrt{\sigma_0^2+ \frac{32 M r_0}{\mu}}}{2} \geq \frac{1}{2}\left(\sigma_0+\sqrt{\frac{M r_0}{\mu}}\right).
	\end{equation}
	Combining Eqs.~\eqref{eq:qm-upper} and \eqref{eq:qm-lower}, we obtain $q_m = \Theta\left(\sigma_0+\sqrt{\frac{M r_0}{\mu}}\right)$.
	Finally, Eq.~\eqref{eq:easyread-convergence} can be derived by choosing $q=q_m$ in Lemma~\ref{lem:base} and Lemma~\ref{lem:super}. \hfill$\square$
\end{proof}

\begin{remark}
	Note that Lemma~\ref{lem:base} and Lemma~\ref{lem:super} require that 
	\[ \frac{M r_0}{\mu}  \leq \frac{q(1-q)}{8} \left(\frac{q}{q+1}-\sigma_0\right) \leq \frac{1}{32} \left(\frac{q}{q+1}-\sigma_0\right), \]
	showing that $\frac{32 M r_0}{\mu} + \sigma_0 \leq \frac{q}{q+1} = O(1)$. Hence we could assume Eq.~\eqref{eq:init-super} in Theorem~\ref{thm:good} holds without loss of generality.
\end{remark}

\subsection{Convergence Rate of Broyden's ``Bad'' Scheme}\label{subsec:bad}

Now we turn to Broyden's ``bad'' scheme. 
Similarly, the Broyden's ``bad'' scheme can be represented by
\begin{equation} \label{eq:h-up}
\left\{
\begin{aligned}
& \vx_{k+1} = \vx_{k} - \mH_{k} \mF(\vx_{k}), \\
& \mH_{k+1} = \mH_k + \frac{\left(\mJ_k^{-1} - \mH_k\right)\vy_k\vy_k^\top}{\vy_k^\top\vy_k}=\text{Broyd}\left(\mH_k, \mJ_k^{-1}, \vy_k\right)
\end{aligned}
\right.
\end{equation}
using $\vu_k\stackrel{\eqref{eq:int-J}}{=} \mJ_k^{-1}\vy_k$ if $\mJ_k$ is nonsingular. 
Now we first provide a local linear convergence rate of $r_k$.
\begin{lemma}\label{lem:base-bad}
	Suppose the objective function $\mF(\cdot)$ satisfies Assumption~\ref{ass:lisp} and sequences $\{\vx_k\}_{k=0}^{+\infty}$ and $\{\mH_k\}_{k=0}^{+\infty}$ generated by Algorithm~\ref{algo:broyden-bad} satisfy Assumption~\ref{ass:optimal}.   
	If there exists some $q\in(0,1)$ such that $r_0$ and $\tau_0$ satisfy
	\begin{equation}\label{eq:init-cond-bad}
		\tau_0 \leq q, \quad  \frac{M R_0}{\mu^2} \leq \frac{q(1-q)}{6} \cdot \left(q-\tau_0\right),
	\end{equation}
	then it holds that for all $k\geq 0$,
	\begin{equation}\label{eq:sigma-bound-bad}
		\tau_k^2 \leq \tau_0^2 + \frac{1+q}{1-q}\left(\frac{2 M R_0}{\mu^2}+\frac{2 M^2 R_0^2}{\mu^4}\right) \leq q^2
	\end{equation}
	and 
	\begin{equation}\label{eq:linear-con-bad}
		r_{k+1} \leq q r_k.
	\end{equation}
\end{lemma}

Next, we provide an explicit superlinear convergence rate as follows.

\begin{lemma}\label{lem:super-bad}
	Suppose the objective function $\mF(\cdot)$ and sequences $\{\vx_k\}_{k=0}^{+\infty}$ and $\{\mH_k\}_{k=0}^{+\infty}$ satisfy the properties described in Lemma~\ref{lem:base-bad}.  
	Let distances $r_0$ and $\tau_0$ satisfy Eq.~\eqref{eq:init-cond-bad} with $q \leq 1/2$.
	Then it holds that for all $k\geq 1$, 
    \[ F_{k} \leq \left[\frac{10 q^2}{k}\right]^{k/2} F_0. \]
    Moreover, we can replace $F_k$ to $f_k$ and conclude that
    \[ f_{k} \leq \varkappa \left[\frac{10 q^2}{k}\right]^{k/2} f_0. \]
\end{lemma}

Similar to the convergence analysis of the ``good'' scheme, we define $q_m$ given the suitable $\tau_0$ and $r_0$ as 
\begin{equation*}
	q_m := \min q, \ s.t.\ \tau_0 \leq q, \ \frac{M R_0}{\mu^2} \leq \frac{q(1-q)}{6} \cdot \left(q-\tau_0\right).
\end{equation*}
Such a $q_m$ needs a pair of proper $\tau_0, r_0$ to guarantee nonempty feasible set, and $q_m$ determines the optimal linear and superlinear convergence rate provided by Lemmas~\ref{lem:base-bad} and~\ref{lem:super-bad}. Then we could obtain the concrete rates expressed by initialization directly.
We leave the proof to Appendix \ref{app:thmbad} because the detail is similar to Theorem \ref{thm:good}. 

\begin{theorem}\label{thm:bad}
	Suppose the objective function $\mF(\cdot)$ satisfies Assumption~\ref{ass:lisp}, and sequences $\{\vx_k\}_{k=0}^{+\infty}$ and $\{\mH_k\}_{k=0}^{+\infty}$ generated by Algorithm~\ref{algo:broyden-bad} satisfy Assumption~\ref{ass:optimal}.  
	If $r_0$ and $\tau_0$ satisfy
	\begin{equation}\label{eq:init-super-bad}
		\frac{24 M R_0}{\mu^2} + \tau_0 \leq \frac{1}{2},
	\end{equation}
   then it holds that $q_m$ is well-defined and $q_m = \Theta\left(\tau_0+\sqrt{\frac{ M R_0}{\mu^2}}\right)$, and we have 
	\begin{equation}\label{eq:easyread-convergence-bad}
		\begin{aligned}
			r_{k+1} &\leq 4\left(\tau_0+\sqrt{\frac{ M R_0}{\mu^2}}\right)r_k, \forall k\geq 0,  \\
		    f_{k} & \leq \varkappa \left[\frac{13\left(\tau_0+\sqrt{\frac{M R_0}{\mu^2 }}\right)}{\sqrt{k}}\right]^{k}f_0, \forall k\geq 1.
		\end{aligned}
	\end{equation}
\end{theorem}

\begin{remark}
	Note that Lemma~\ref{lem:base-bad} and Lemma~\ref{lem:super-bad} require that 
	\[ \frac{M R_0}{\mu^2} \leq \frac{q(1-q)}{6} \cdot \left(q - \tau_0\right) \leq \frac{1}{24} \left(q-\tau_0\right), \]
	showing that $\frac{24 M R_0}{\mu^2} + \tau_0 \leq q = O(\frac{1}{\varkappa})$. Hence we could assume Eq.~\eqref{eq:init-super-bad} in Theorem~\ref{thm:bad} holds without loss of generality.
\end{remark}

\subsection{Comparison Between Broyden's ``Good'' and ``Bad'' Schemes}\label{subsec:compare}

Theorems~\ref{thm:good} and~\ref{thm:bad} provide the conditions to guarantee the convergence and the explicit superlinear convergence rates of Broyden's ``good'' and ``bad'' schemes.
We now give a detailed comparison between the ``good'' and ``bad'' schemes and attempt to present some useful insights to explain the differences of these two schemes.

First, Theorems~\ref{thm:good} and~\ref{thm:bad} show that the local regions to achieve superlinear convergence rates of the ``good'' and ``bad'' schemes are quite similar.
The ``good'' scheme requires that $\frac{32M r_0}{\mu} + \sigma_0 \leq \frac{1}{3}$ which implies that $\frac{Mr_0}{\mu} = O(1)$ and $\sigma_0 = O(1)$.
In contrast, the local region for the ``bad'' scheme needs $\frac{24 M R_0}{\mu^2} + \tau_0 \leq \frac{1}{2}$, implying that $\frac{MR_0}{\mu^2} = O(1)$, and $\tau_0 = O(1)$. 
Assume the initial Jacobian can be represented as $ \mB_0 = s \mJ_* $ with $0 < s < +\infty$.
By the definitions of $\sigma_0$ and $\tau_0$, we can obtain $\sigma_0 = \sqrt{n} \cdot  \lvert s-1 \rvert $ and $\tau_0 = \sqrt{n} \cdot \lvert s^{-1}-1\rvert$. 
Thus, the local region is $\lvert s-1 \rvert = O(\frac{1}{\sqrt{n}})$ for both ``good'' and ``bad'' schemes.

\begin{wrapfigure}{r}{0.3\textwidth}
	\centering
	\vspace{-15pt}
	\includegraphics[width=0.2\textwidth]{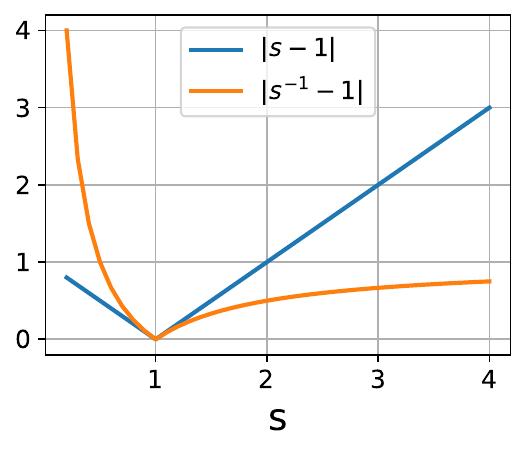}
	\vspace{-10pt}
	\caption{$\lvert s-1 \rvert$ vs. $\lvert \frac{1}{s}-1 \rvert$.}
	\label{fig:s-plot}
	\vspace{-40pt}
\end{wrapfigure}

\begin{table*}[t]
	\centering
	\caption{Comparison of the two schemes with the assumption $\mB_0 = s \mJ_*, s>0$.}
	\begin{tabular}{ccc}
		\toprule
		$\mB_0=s\mJ_*, s>0$ & The ``good'' scheme & The ``bad'' scheme \\
		\midrule
		Local Region & $\frac{Mr_0}{\mu} = O(1), \lvert s -1 \rvert=O(\frac{1}{\sqrt{n}})$ & $\frac{MR_0}{\mu^2} = O(1), \lvert s-1 \rvert =O(\frac{1}{\sqrt{n}})$ \\  
		\midrule
		Superlinear Rate & $\left(\dfrac{O\left(\sqrt{n}\lvert s -1 \rvert+\sqrt{Mr_0/\mu}\right)}{\sqrt{k}}\right)^{k}$ & $\left(\dfrac{O\left(\sqrt{n}\lvert s^{-1} -1 \rvert+\sqrt{MR_0 / \mu^2}\right)}{\sqrt{k}}\right)^{k}$ \\
		\midrule
		Suitable Scene & $s \leq 1, \kappa \gg 1$ & $s \geq 1, \kappa \approx 1$ \\
		\bottomrule
	\end{tabular}
	\label{table:res}
\end{table*}

Second, Theorems~\ref{thm:good} and~\ref{thm:bad} show that Broyden's ``good'' and ``bad'' schemes have the local superlinear convergence rates $\big(\frac{O(\sigma_0+\sqrt{M r_0 / \mu})}{\sqrt{k}}\big)^{k}$ and $\big(\frac{O(\tau_0+\sqrt{M R_0/ \mu^2})}{\sqrt{k}}\big)^{k}$, respectively.
We can observe that the initial distances $r_0, \sigma_0$ (or $R_0, \tau_0$) determine the convergence speeds of the ``good'' and ``bad'' schemes.

Now we discuss how the initial error of the approximate Jacobian leads to difference performance of the ``good'' and ``bad'' schemes.
Assume that the distance between the initial Jacobian and the exact one dominates the convergence rates and $ \mB_0 = s \mJ_* $ with $0 < s < +\infty$, which implies that $\sigma_0 = \sqrt{n} \cdot \lvert s-1 \rvert $, $\tau_0 = \sqrt{n} \cdot \lvert s^{-1}-1\rvert$, and the explicit rates of the ``good'' and ``bad'' schemes reduce to $\left(\frac{O(\sqrt{n}\lvert s-1\rvert)}{\sqrt{k}}\right)^{k}$ and $\left(\frac{O(\sqrt{n}\lvert s^{-1} - 1\rvert)}{\sqrt{k}}\right)^{k}$. 
Thus, it holds that $\sigma_0 \le \tau_0$ when $s\le 1$, which implies the ``good'' scheme outperforms the ``bad'' scheme.
On the other hand, when $s>1$, it holds that $\tau_0 < \sigma_0$ which shows that the ``bad'' scheme achieves  better performance.
In fact, these phenomena have been conjectured in the existing work \cite[Section 3.1]{martinez2000practical}. 
Although the ``bad'' scheme may achieve faster convergence rate than the ``good'' scheme when $s>1$, its convergence properties will deteriorate severely when  $s$ is close to zero. 
As shown in Figure \ref{fig:s-plot}, $\lvert s^{-1}-1 \rvert$ is highly nonlinear when $0<s<1$.
In contrast, $\lvert s-1 \rvert$ is piece-wise linear.
Thus, the convergence rate of the ``good'' is much robust to the value of $s$.
The above discussion may provide some explanation of why the ``good'' scheme is more robust than the ``bad'' scheme in practice.

We also try to figure out how the initial distances $r_0$ and $R_0$ affect the convergence rates of the ``good'' and ``bad'' schemes.
We assume that the initial approximate Jacobian $\mB_0$ is sufficiently close to $\mJ_{*}$, that is, the values of $\sigma_0$ and $\tau_0$ are close to zeros.
In this case, noting that $\mu r_0 \leq R_0 \leq L r_0$ by Eq.~\eqref{eq:rR}, we get $Mr_0/\mu \leq MR_0/\mu^2 \leq \kappa \cdot Mr_0 /\mu$. 
Thus the ``good'' scheme always achieve a faster convergence rate than the ``bad'' scheme, particularly for a large $\kappa$.
Indeed, the classical BFGS and DFP methods share a similar idea as Broyden's ``good'' and ``bad'' methods, which update the Hessian matrix and inverse Hessian directly. The current superlinear rates \cite{rodomanov2021new,rodomanov2021rates} also reveal that BFGS method is faster than DFP method (with a condition number $\kappa$).

Overall, we summarize the comparison discussed above to Table~\ref{table:res} based on the assumption that $\mB_0 = s \mJ_*$ with $s>0$.
\section{Experiments}\label{sec:exp}

In this section we mainly validate our theoretical findings on the initial approximate Jacobian shown in Table \ref{table:res}. 
To monitor different suitable initial points, we may use Newton's method to obtain a (nondegenerate) solution $\vx_*$ and then choose an initial $\vx_0$ around $\vx_*$ with a tiny perturbation.   
We employ the initial approximate Jacobian $\mB_0 = s\mJ(\vx_*)$ with $s>0$, or $\mB_0 = s\mJ(\vx_0)$ which is more common and practical in real problems \cite{nocedal2006numerical}. 

\begin{figure}[t]
	\centering
	\hspace{-10pt}
	\begin{subfigure}[b]{0.5\textwidth}
	\includegraphics[width=\linewidth]{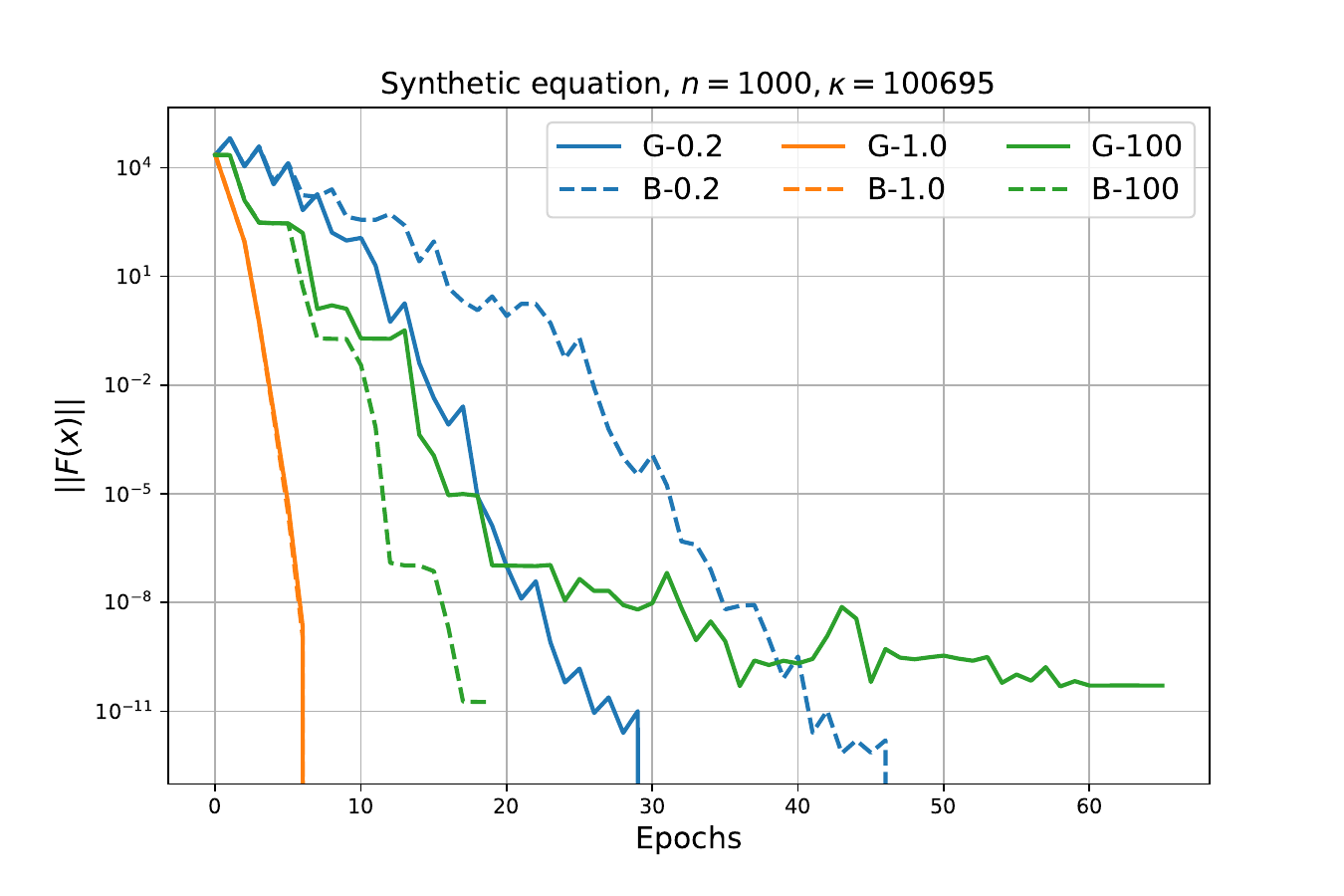}
	\caption{Synthetic equation with $\mB_0=s\nabla^2 f(\vx_*)$.} \label{fig:broyden-syn-equation}
	\end{subfigure}
	\hspace{-20pt}
	\begin{subfigure}[b]{0.5\textwidth}
	\includegraphics[width=\linewidth]{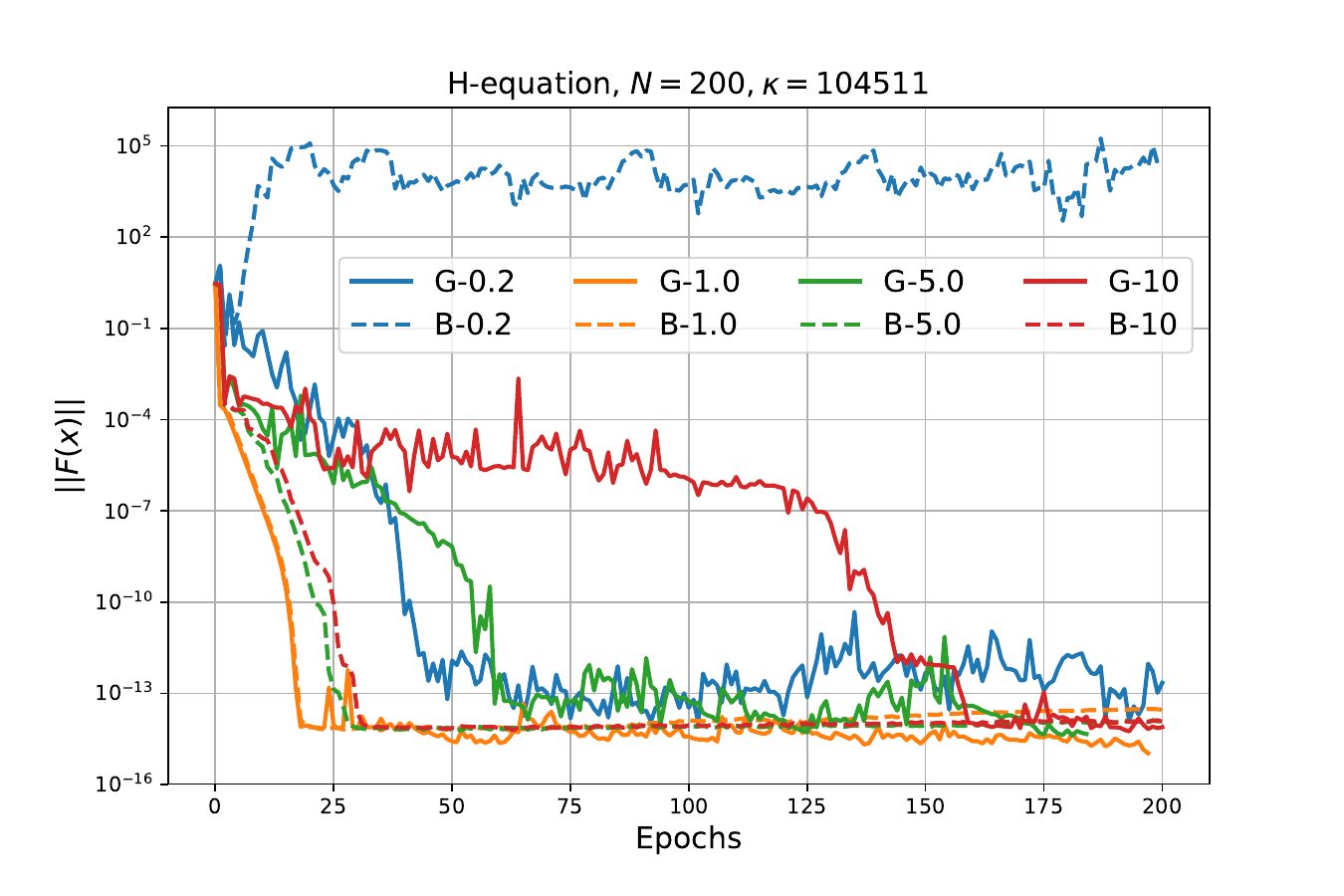}
	\caption{H-equation with $\mB_0=s\nabla^2 f(\vx_0)$.} \label{fig:broyden-H-equation}
	\end{subfigure}
	\caption{Comparison of Broyden's ``good'' and ``bad'' schemes for the (a) Synthetic equation and (b) Chandrasekhar H-equation. We choose $ \mB_0=s\nabla^2 f(\vx_*)$ or $s\nabla^2 f(\vx_0)$ with various $s>0$ shown in legends. Here the legend ``G(B)-0.2'' means using good (bad) method with $s=0.2$.
	The constant in H-equation is $c=1{-}10^{-10}$. 
	We also list the dimension $n$ or $N$ and the condition number $\varkappa$ in each title.}\label{exp-figs}
\end{figure}

\textbf{Synthetic nonlinear equation.}
We first solve the following simple synthetic nonlinear equation:
\[ \mF(\vx) = \mA (\vx \odot \vx-\bm{1}_n)= \bm{0}, \ \vx \in \sR^n. \]
where $\mA\in \sR^{n \times n}$ is a predefined non-singular matrix, $\odot$ is the element-wise product and $\bm{1}_n$ is the $n$-dimensional vector filled with ones. Clearly, a nondegenerate solution is $\vx_*=\bm{1}_n$, because the Jacobian $\mJ(\vx_*) = 2\mA$. 
Assumption~\ref{ass:lisp} is satisfied because $\mJ(\vx) = 2\mA \; \mathrm{diag}\{\vx\}$ where $\mathrm{diag}\{\vx\}$ is the diagonal matrix with the elements of vector $\vx$ on the main diagonal.
We use Broyden's method to solve problem $\mF(\vx) = \bm{0}$. 
We choose the shared starting point $\vx_0 = \vx_*+0.1\cdot\norm{\vx_*}\cdot\bm{\epsilon}$, where $\bm{\epsilon} \sim \text{Unif}\left(\gS^{n-1}\right)$ and $\mB_0=s \mJ(\vx_*)$ with various $s>0$. 
We report numerical results in Figure~\ref{fig:broyden-syn-equation}.

Figure~\ref{fig:broyden-syn-equation} shows that the two algorithms perform similarly when $s=1$. But when $s$ is smaller than $1$, the ``good'' scheme outperforms the ``bad'' scheme ($s=0.2$) and the reverse finding appears when $s$ is far away from $1$ ($s=100$). Surprisingly, the ``good'' scheme does not fail even when $s=100$, but it indeed converges much slower.

\textbf{The Chandrasekhar H-equation.}
We also attempt to use Broyden's methods to solve the Chandrasekhar H-equation \cite[Chapter 5.6]{kelley1995iterative}. 
The problem needs to solve the following functional problem:
\[ \mF(\gH)(\mu) = \gH(\mu)-\left(1-\frac{c}{2}\int_{0}^1 \frac{\mu \gH(\nu)d\nu}{\mu+\nu}\right)^{-1} = 0, \]
where $\gH: [0, 1]\to \sR$ and $c\in(0, 1)$. We will discretize the equation with the composite midpoint rule. Here we
approximate integrals on $[0, 1]$ by
\[ \int_{0}^1 f(\mu)d\mu \approx \frac{1}{N}\sum_{j=1}^N f(\mu_j), \]
which $\mu_i = (i-0.5)/N, i\in[N]$. The resulting discrete problem is
\[ F_i(\vx) = x_i - \left(1-\frac{c}{2N}\sum_{j=1}^N\frac{\mu_ix_j}{\mu_i+\mu_j}\right)^{-1} = 0, i \in [N], \]
where $\vx = (x_1, \dots, x_N)^\top \in \sR^N, \mF(\vx) = (F_1(\vx), \dots, F_N(\vx))^\top \in \sR^N$.
The Jacobian of $\mF(\vx)$ is $\mJ(\vx) = [J_{i k}(\vx)] \in\sR^{N \times N}$ with
\[ J_{i k}(\vx) = \frac{\partial F_i(\vx)}{\partial x_k} = \mathbbm{1}_{\{i=k\}}-\frac{c}{2N} \cdot \frac{\mu_i}{\mu_i+\mu_k} \cdot \left(1-\frac{c}{2N}\sum_{j=1}^N\frac{\mu_ix_j}{\mu_i+\mu_j}\right)^{-2}, \]
where $\mathbbm{1}_{\{i=k\}}$ is the indicator function that equals to $1$ only when $i$ equals to $k$.
It is known \cite{mullikin1968some} that the discrete analog has
solutions for $ c \in (0, 1) $. Thus it is differentiable, and locally Lipschitz continuous restricting in a bounded domain near the optimal solution $\vx_*$, which satisfies Assumption~\ref{ass:lisp}.
We also adopt Broyden's method to solve problem $\mF(\vx) = \bm{0}$ with the shared starting point $\vx_0 = \vx_*+0.1\cdot\norm{\vx_*}\cdot\bm{\epsilon}$, where $\bm{\epsilon} \sim \text{Unif}\left(\gS^{N-1}\right)$ and  $\mB_0=s \mJ(\vx_0)$ with various $s>0$. 
We report numerical results in Figure~\ref{fig:broyden-H-equation}. 

Figure~\ref{fig:broyden-H-equation} also shows the different performance between ``good'' and ``bad'' schemes when $s<1$ or $s>1$, while the two algorithms perform similarly when $s=1$.
Moreover, the ``bad'' scheme fails when $s=0.2$, showing that the ``bad'' scheme is nonrobust to underestimated Jacobian.
However, the ``good'' scheme does not fail even for a large $s$ as well.

In summary, by the numerical results shown in Figure~\ref{exp-figs}, we may claim that the ``good'' and ``bad'' schemes indeed have their own applicable scenes. 
The rough conclusion is that ``good'' scheme is more suitable for underestimated initial approximate Jacobian (i.e., $s<1$), while the ``bad'' scheme is applied to overestimated initial approximate Jacobian (i.e., $s\geq1$).
On the other hand, when we do not know the initial approximate Jacobian is underestimated or overestimated, the ``good'' scheme is preferred because its performance is robust to the choice of the initial approximate Jacobian (Fig.~\ref{fig:s-plot}).

\section{Comparison with Existing Explicit Rates of Quasi-Newton Methods}\label{sec:compare}
Since there have been many recent works \cite{rodomanov2021greedy,rodomanov2021rates,rodomanov2021new,jin2020non,lin2021faster,ye2021explicit} on the explicit non-asymptotic superlinear rates of quasi-Newton methods, we list our difference in detail. 
One important observation is that previous work mainly focused on unconstrained optimization rather than solving the nonlinear equations directly.
Despite the many similarities between nonlinear equations and unconstrained optimization algorithms, there are also some important differences.
\begin{enumerate}
    \item First, the Jacobian in solving a nonlinear equation system is commonly asymmetric while the Hessian in the unconstrained optimization is symmetric.
    And we can convert unconstrained convex optimization $\min_{\vx} f(\vx)$ to solving a nonlinear equation $\nabla f(\vx) =\bm{0}$.
    Thus, quasi-Newton methods for solving nonlinear equations have a wider range of applications than for unconstrained optimization.
    \item Second, some methods would define a merit function \cite[Chapter 11.2]{nocedal2006numerical},
    e.g., $\norm{\mF(\vx)}_2^2, \norm{\mF(\vx)}_1$, to convert solving the nonlinear equations to an unconstrained optimization problem. 
    Whereas, various merit functions can be used in nonlinear equations, and choosing a proper one is also a challenge. 
    Moreover, all of which may have some drawbacks \cite[Section 11.0]{nocedal2006numerical}. 
    \begin{enumerate}
        \item For example, we need to know the expression of $\mJ(\vx)$ to obtain the gradient of a merit function generally, but we only need the function value $\mF(\vx)$ for solving nonlinear equations in Algorithms \ref{algo:broyden} and \ref{algo:broyden-bad}.
        \item Additionally, the common used merit function, $\norm{\mF(\vx)}_2^2$, may suffer from the numerical unstability. 
        Particularly, to achieve $\norm{\mF(\vx)}_2 \leq \epsilon$ for some small $\epsilon>0$, Broyden's methods could directly solve the problem. However, it needs to solve the problem $\norm{\mF(\vx)}_2^2 \leq \epsilon^2$ by adopting the merit function, which may be unsolvable due to the high precision requirement ($\epsilon^2$) of the solution.
        Moreover, such merit function would bring a larger condition number ($\varkappa^2$) than solving the original nonlinear equations directly\footnote{$\nabla^2 \left(\norm{\mF(\vx_*)}_2^2/2\right) = \mJ_*^\top \mJ_* + \sum_{i=1}^n F_i(\vx_*)\nabla^2 F_i(\vx_*)= \mJ_*^\top\mJ_*$ by $F_i(\vx_*)=0, \forall i \in [n]$.}.
    \end{enumerate}
\end{enumerate}

Furthermore, the drawback of our analysis is that we needs an additional condition on the initial approximate Jacobian matrix following \cite[Theorem 8.2.2]{dennis1996numerical} and \cite[Theorem 11.5]{nocedal2006numerical}.
Such a requirement is unavoidable because Broyden's update does not preserve the non-singularity of approximate Jacobian matrix, but classical quasi-Newton methods (such as BFGS and DFP) indeed preserve the positive definiteness (e.g., \cite[Lemma 2.1]{rodomanov2021rates}).
Though there are several variations of original Broyden's methods (\cite{martinez2000practical,dennis1977quasi}) to overcome such a drawback of non-singularity, the explicit rates maybe loose. Thus, we still focus on the original version of Broyden's methods providing that the initial approximate Jacobian is good.

Due to differences mentioned above, we consider that our work is different from existing works on explicit superlinear convergence rates of quasi-Newton methods 
though our work adopts the semblable ideas of \cite{jin2020non,lin2021faster,rodomanov2021greedy,rodomanov2021new,rodomanov2021rates,ye2021explicit}.

\section{Conclusion} \label{sec:conclusion}
In this work, we have studied the behavior of the famous Broyden's method for solving the general nonlinear equations.
We have presented explicit local superlinear convergence rates of Broyden's method. 
Our theoretical results of Broyden's method, including Broyden's ``good'' and ``bad'' schemes, reveal the factors that determine the convergence rates of these two schemes and provide some important insights on the performance difference between these two versions.
Specifically, our work shows that Broyden's ``good'' scheme has a wide scope of suitable initial approximate Jacobian, while the ``bad'' scheme suffers from underestimated approximate Jacobian. 
However, once the initialization is close enough to a solution, these two schemes could perform similarly with comparable superlinear convergence rates. 
We hope our understanding could give a more quantitative understanding of Broyden's method, as well as the family of quasi-Newton methods.

\section*{Acknowledgments}

The authors would like to thank Chris Junchi Li for helpful discussion and pointing out the relation with Broyden's ``bad'' scheme on the paper.

\appendix


\section{Auxiliary Lemmas}
\begin{lemma}[Part of \cite{dennis1996numerical} Theorem 3.1.3]
Let $\mA, \mB \in\sR^{n\times n}$. Then 
\[ \norm{\mA\mB}_F \leq \norm{\mA} \cdot \norm{\mB}_F. \]
\end{lemma}

\begin{lemma}[Extension of \cite{dennis1996numerical} Theorem 3.1.4] \label{aux-lemma2}
Let $\|\cdot \|$ be any norm defined on $\sR^{n\times n}$ that obeys
$ \|\mA\mB\| \leq \|\mA\|\cdot \|\mB\| $ for any $\mA,\mB \in\sR^{n\times n}$
and $\|\mI_n\|=1$, and let $\mE \in\sR^{n\times n}$. 
If $\norm{\mE} < 1$, then $\left(\mI_n-\mE\right)^{-1}$ exists with
\begin{equation}\label{eq:inv-i-e}
	\norm{\left(\mI_n-\mE\right)^{-1}} \leq \frac{1}{1-\norm{\mE}},
\end{equation}
and
\begin{equation}\label{eq:inv-i-e-i}
	\norm{\left(\mI_n-\mE\right)^{-1}-\mI_n} \leq \frac{\norm{\mE}}{1-\norm{\mE}}.
\end{equation}
\end{lemma}

\begin{proof}
Since $\norm{\mE} < 1$, we have $0 \leq \liminf_{k\to \infty}\norm{\mE^k} \leq \limsup_{k\to \infty}\norm{\mE^k} \leq \lim_{k\to \infty}\norm{\mE}^k =0 $, thus $\lim_{k\to \infty}\norm{\mE^k}=0$ and $\lim_{k\to \infty} \mE^k = \bm{0}$, leading to $\lim_{k\to \infty} \left(\mI_n-\mE\right)\left(\sum_{i=0}^k \mE^i\right) = \lim_{k\to \infty} \mI_n-\mE^{k+1} = \mI_n$. Hence
\[ \left(\mI_n-\mE\right)^{-1} = \sum_{i=0}^\infty \mE^i = \mI_n+\mE+\mE^2+\dots. \]
Therefore $\left(\mI_n-\mE\right)^{-1}$ exists, and
\[ \norm{\left(\mI_n-\mE\right)^{-1}} = \norm{\sum_{i=0}^\infty \mE^i} \leq \sum_{i=0}^\infty \norm{ \mE^i} \leq \sum_{i=0}^\infty \norm{ \mE}^i = \frac{1}{1-\norm{\mE}}, \]
and
\[ \norm{\left(\mI_n-\mE\right)^{-1}-\mI_n} = \norm{\sum_{i=1}^\infty \mE^i} \leq \sum_{i=1}^\infty \norm{\mE^i}  \leq \sum_{i=1}^\infty \norm{ \mE}^i = \frac{\norm{\mE}}{1-\norm{\mE}}. \]
The proof is finished. \hfill$\square$
\end{proof}

Additionally, we give the well-known estimator of the objective function and its Jacobian around $\vx_*$ under Assumption~\ref{ass:lisp}.
\begin{lemma}[Analogy to \cite{dennis1996numerical} Lemmas 4.1.15 and 4.1.12]\label{lemma:prop}
Let $\mF(\vx)$ satisfy Assumption~\ref{ass:lisp}, and  $\mJ := \int_{0}^{1} \mJ\left(\vx+t\vu\right)dt$ for some $\vx, \vu \in\sR^n$. Then the following inequalities hold:
\begin{align}
	& \norm{\mJ - \mJ(\vx_*)} \leq \frac{M}{2} \left(\norm{ \vx- \vx_*}+\norm{ \vx- \vx_* + \vu}\right), \label{eq:lemma-j1} \\
	& \|\mF(\vx) - \mJ(\vx_*)(\vx-\vx_*) \| \leq \frac{M}{2} \norm{ \vx-\vx_* }^2. \label{eq:lemma-j2}
\end{align}		
\end{lemma}

\begin{proof}
It holds that
\begin{eqnarray*}
	&&\|\mJ - \mJ(\vx_*)\| \leq \| \int_{0}^{1} \mJ\left(\vx + t\vu\right) - \mJ(\vx_*)dt \| \leq \int_{0}^{1} \norm{\mJ\left(\vx + t\vu\right) - \mJ(\vx_*)} dt \\ \nonumber
	&\stackrel{\eqref{eq:lisp}}{\leq}&  M \int_{0}^{1} \norm{ \vx + t\vu - \vx_*} d t \leq M  \int_{0}^{1} \left[(1-t)\norm{ \vx-\vx_*} + t\norm{ \vx - \vx_* + \vu} \right] dt \\ \nonumber
	&=& \frac{M}{2} \left(\norm{ \vx- \vx_*}+\norm{ \vx- \vx_* + \vu}\right),
\end{eqnarray*}
where the fourth inequality uses the convexity of $\norm{\cdot}$.
Thus, Eq.~\eqref{eq:lemma-j1} is proved. \\
Next, we denote $\vs = \vx-\vx_*$. Then from Assumption~\ref{ass:lisp}, we get
\begin{eqnarray*}
	&&\lefteqn{ \norm{\mF(\vx) - \mJ(\vx_*)(\vx-\vx_*)} = \norm{\mF(\vx) - \mF(\vx_*) - \mJ(\vx_*)(\vx-\vx_*)} } \\ \nonumber
	&=& \| \int_{0}^{1} \left[\mJ\left(\vx_*+t\vs\right) - \mJ(\vx_*)\right]\vs dt \| \leq \int_{0}^{1} \norm{\mJ\left(\vx_*+t\vs\right) - \mJ(\vx_*)} \cdot \norm{\vs} dt \\ \nonumber
	&\stackrel{\eqref{eq:lisp}}{\leq} & \int_{0}^{1} Mt \norm{\vs}^2 dt = \frac{M}{2} \norm{ \vx-\vx_*}^2,
\end{eqnarray*}
which shows that Eq.~\eqref{eq:lemma-j2} holds.  \hfill$\square$
\end{proof}

\section{Missing Proofs in Section~\ref{subsec:good}}

Before diving into proving Lemmas~\ref{lem:base} and~\ref{lem:super}, we are going to give two important lemmas which are the foundation to proofs of  Lemmas~\ref{lem:base} and~\ref{lem:super}.
First, we build up the relation between $\sigma_{k+1}$ and $\sigma_k$  below.

\begin{lemma}\label{lemma:sigma-up}
Following Algorithm~\ref{algo:broyden} and denoting $\overline{\mJ_k} := \mJ_*^{-1}(\mJ_k - \mJ_*)$, then we have for all $k \geq 0$,
\begin{equation}\label{eq:sigma-update1}
	\sigma_{k+1}^2 \leq \sigma_k^2 -  \frac{\norm{\mJ_*^{-1}\left(\mJ_k-\mB_k\right)\vu_k}^2}{\norm{\vu_k}^2} + 2 \left(\sigma_k\norm{\overline{\mJ_k}} + \norm{\overline{\mJ_k}}^2\right).
\end{equation}
Moreover, we can introduce $f_k$ into the above inequality:
\begin{equation}\label{eq:sigma-update2}
	\begin{aligned}
		\sigma_{k+1}^2 \leq \sigma_k^2 - \frac{f_{k+1}^2}{f_k^2} \cdot \frac{\norm{\mJ_*^{-1}\mB_k\vu_k}^2}{\norm{\vu_k}^2} + 2 \left(\sigma_k\norm{\overline{\mJ_k}} + \norm{\overline{\mJ_k}}^2\right).
	\end{aligned}
\end{equation}
\end{lemma}

\begin{proof}
Following the definition of $\sigma_{k+1}$, we have 
\begin{equation}\label{eq:bde}
	\overline{\mB}_{k+1} \stackrel{\eqref{eq:jac-measure}}{=} \mJ_*^{-1}\left[\mB_k {+}  \frac{\left(\mJ_k {-} \mB_k\right)\vu_k\vu_k^\top}{\vu_k^\top\vu_k} {-} \mJ_{*}\right] = \overline{\mB}_k +  \frac{\mJ_*^{-1}\left(\mJ_k {-} \mB_k\right)\vu_k\vu_k^\top}{\vu_k^\top\vu_k}.
\end{equation}
Taking the Frobenius norm in both  sides of the above equation, we obtain
\begin{eqnarray}\label{eq:part1}
	&& \sigma_{k+1}^2 \stackrel{\eqref{eq:jac-measure}}{=} \norm{\overline{\mB}_{k+1}}_F^2 \stackrel{\eqref{eq:bde}}{=} \norm{\overline{\mB}_k +  \frac{\mJ_*^{-1}\left(\mJ_k-\mB_k\right)\vu_k\vu_k^\top}{\vu_k^\top\vu_k}}_F^2 \\ \nonumber
	&=& \norm{\overline{\mB}_k}_F^2 + \frac{\norm{\mJ_*^{-1} \left(\mJ_k-\mB_k\right)\vu_k}^2}{\norm{\vu_k}^2} + 2 \cdot \frac{\vu_k^\top \left(\mB_k -\mJ_{*}\right)^\top \mJ_*^{-\top}\mJ_*^{-1} \left(\mJ_k-\mB_k\right)\vu_k}{\vu_k^\top\vu_k} \\ \nonumber
	&\stackrel{\eqref{eq:jac-measure}}{=} & \sigma_k^2 + \frac{\norm{\mJ_*^{-1} \left(\mJ_k-\mB_k\right)\vu_k}^2}{\norm{\vu_k}^2} + 2 \cdot \frac{\vu_k^\top \left(\mB_k-\mJ_k\right)^\top \mJ_*^{-\top}\mJ_*^{-1} \left(\mJ_k-\mB_k\right)\vu_k}{\vu_k^\top\vu_k} \\ \nonumber
	&&+ 2 \cdot \frac{\vu_k^\top \left(\mJ_k-\mJ_{*}\right)^\top \mJ_*^{-\top}\mJ_*^{-1} \left(\mJ_k-\mB_k\right)\vu_k}{\vu_k^\top\vu_k} \\ \nonumber
	&=& \sigma_k^2 - \frac{\norm{ \mJ_*^{-1}\left(\mJ_k-\mB_k\right)\vu_k}^2}{\norm{\vu_k}^2} + 2 \cdot \frac{\vu_k^\top \left(\mJ_k-\mJ_{*}\right)^\top \mJ_*^{-\top}\mJ_*^{-1} \left(\mJ_k-\mB_k\right)\vu_k}{\vu_k^\top\vu_k}.
\end{eqnarray}
We also have the last term in the above
\begin{eqnarray}\nonumber
	&&\frac{\vu_k^\top \left(\mJ_k-\mJ_{*}\right)^\top \mJ_*^{-\top}\mJ_*^{-1} \left(\mJ_k-\mB_k\right)\vu_k}{\vu_k^\top\vu_k} \leq \frac{\norm{ \mJ_*^{-1}\left(\mJ_k - \mJ_{*}\right) \vu_k}}{\|\vu_k\|}\cdot \frac{\norm{\mJ_*^{-1}\left(\mJ_k-\mB_k\right)\vu_k}}{\|\vu_k\|} \\ \nonumber
	&\leq& \norm{\mJ_*^{-1} \left(\mJ_k - \mJ_{*}\right) } \cdot \left[ \norm{\mJ_*^{-1}\left(\mB_k-\mJ_{*}\right)} + \norm{\mJ_*^{-1}\left(\mJ_{*}-\mJ_k\right)}\right] \stackrel{\eqref{eq:jac-measure}}{\leq} \sigma_k\norm{\overline{\mJ_k}} + \norm{\overline{\mJ_k}}^2. \label{eq:part2}
\end{eqnarray}
Combining Eqs.~\eqref{eq:part1} and \eqref{eq:part2}, we can obtain Eq.~\eqref{eq:sigma-update1}.

Finally, from Steps 3 and 4 in Algorithm \ref{algo:broyden}, we have
\begin{equation}\label{eq:fjb}
    \mB_k\vu_k = -\mF(\vx_k), \ \mF(\vx_{k+1}) = \vy_k+\mF(\vx_k) \stackrel{\eqref{eq:int-J}}{=} \left(\mJ_k-\mB_k\right) \vu_k.
\end{equation}
Hence, we obtain 
\begin{equation}\label{eq:theta-lam}
    \frac{f_{k+1}} {f_k} \stackrel{\eqref{eq:lam}}{=} \frac{\norm{\mJ_*^{-1}\mF(\vx_{k+1})}}{\norm{\mJ_*^{-1}\mF(\vx_{k})}} \stackrel{\eqref{eq:fjb}}{=} \frac{\norm{\mJ_*^{-1}\left(\mJ_k-\mB_k\right)\vu_k}}{\norm{\mJ_*^{-1}\mB_k\vu_k}}.
\end{equation}
Replacing Eq.~\eqref{eq:theta-lam} to Eq.~\eqref{eq:sigma-update1}, we can obtain the results of Eq.~\eqref{eq:sigma-update2}. \hfill$\square$
\end{proof}


The remaining term $\norm{\overline{\mJ_k}}$ in Lemma~\ref{lemma:sigma-up} can be bounded by the distance $\norm{\vx_k-\vx_*} + \norm{\vx_{k+1}-\vx_*}$ from Lemma~\ref{lemma:prop}. Hence, we also need the property of the residual $r_k$.
\begin{lemma}\label{lemma:r-up}
	Following Algorithm~\ref{algo:broyden}, then we have for all $k \geq 0$,
	\begin{equation}\label{eq:r-update}
		r_{k+1} \leq \frac{\frac{1}{2\mu} \cdot M r_k + \sigma_k}{1-\sigma_k} \cdot r_k, \text{ if } \sigma_k<1,
	\end{equation}	
\end{lemma}
\begin{proof}
	Note that
	\begin{eqnarray*}
		\vx_{k+1}-\vx_* &\stackrel{\eqref{eq:x-up}}{=}&\vx_{k}-\mB_k^{-1}\mF(\vx_k)-\vx_* = \mB_k^{-1}\bigg[\mB_k \left(\vx_k-\vx_*\right)-\mF(\vx_k)\bigg] \\
		&=& \mB_k^{-1}\mJ_*\bigg[\mJ_*^{-1}\left(\mB_k-\mJ_{*}\right) \left(\vx_k-\vx_*\right)-\mJ_*^{-1}\left[\mF(\vx_k)-\mJ_{*}\left(\vx_k-\vx_*\right)\right]\bigg].
	\end{eqnarray*}	
	Taking the spectral norm for the both sides of the above equation, we have that
	\begin{eqnarray}\label{eq:rk-up1} 
		r_{k+1} 
		&\leq& \norm{\mB_k^{-1}\mJ_*} \cdot \bigg[ \norm{\mJ_*^{-1}\left(\mB_k-\mJ_{*}\right)} \cdot r_k + \norm{\mJ^{-1}_*} \cdot \norm{\mF(\vx_k)-\mJ_{*} \left(\vx_k-\vx_*\right)} \bigg] \\ 
		&\stackrel{\eqref{eq:jac-measure}}{\leq}& \norm{\mB_k^{-1}\mJ_*} \cdot \bigg[\sigma_k r_k +\norm{\mJ^{-1}_*} \cdot \norm{\mF(\vx_k)-\mJ_{*} \left(\vx_k-\vx_*\right)} \bigg] \stackrel{\eqref{eq:lemma-j2}}{\leq} \norm{\mB_k^{-1}\mJ_*} \cdot \left[\sigma_k r_k +\frac{M r_k^2}{2\mu} \right].\nonumber
	\end{eqnarray}
	If $\sigma_k<1$, we obtain
	\[ \norm{\mI_n-\mJ_{*}^{-1}\mB_k} = \norm{\mJ_{*}^{-1}\left(\mB_k-\mJ_{*}\right)} \stackrel{\eqref{eq:jac-measure}}{\leq} \sigma_k <1. \]
	Hence from Lemma~\ref{aux-lemma2}, we obtain that $\mJ_{*}^{-1}\mB_k$ is nonsingular, and 
	\begin{equation}\label{eq:inv}
		\norm{\mB_k^{-1}\mJ_*} \stackrel{\eqref{eq:inv-i-e}}{\leq} \frac{1}{1-\norm{\mJ_{*}^{-1}\left(\mB_k-\mJ_*\right)}} \leq \frac{1}{1-\sigma_k}, \; \text{ if } \sigma_k<1.
	\end{equation}
	Finally, we derive Eq.~\eqref{eq:r-update} from Eqs.~\eqref{eq:rk-up1} and \eqref{eq:inv}. \hfill$\square$
\end{proof}

Lemma~\ref{lemma:r-up} shows that  $\{r_k\}$ converges linearly if the initial values $r_0$ and $\sigma_0$ are small enough. 
Hence we could obtain a local linear convergence rate if the starting point $\vx_0$ is close enough to the optimal $\vx_*$ and $\mB_0$ is close enough to $\mJ_{*}$.

\subsection{Proof of Lemma~\ref{lem:base}}
Next, we provide the proof of Lemma~\ref{lem:base} as follows.

\begin{proof}
	Denoting $\overline{\mJ_k}=\mJ_*^{-1}\left(\mJ_k-\mJ_*\right)$ following Lemma \ref{lemma:sigma-up}, we have 
	\begin{equation}\label{eq:over-j}
		\norm{\overline{\mJ_k}} \leq \norm{\mJ_*^{-1}} \cdot \norm{\mJ_k-\mJ_*} \stackrel{\eqref{eq:lemma-j1}}{\leq} \frac{M}{2\mu} \left(\norm{\vx_k-\vx_*}+\norm{\vx_{k+1}-\vx_*}\right) = \frac{M \left(r_k+r_{k+1}\right)}{2\mu}.
	\end{equation}
	Combining with Lemma~\ref{lemma:sigma-up}, we can derive that
	\begin{equation}\label{eq:sk-up}
	    \sigma_{k+1}^2 \stackrel{\eqref{eq:sigma-update1}}{\leq} \sigma_k^2 + 2 \left(\sigma_k\norm{\overline{\mJ_k}} + \norm{\overline{\mJ_k}}^2\right) \stackrel{\eqref{eq:over-j}}{\leq} \sigma_k^2 + 2 \left(\frac{\sigma_k M\left(r_k + r_{k+1}\right)}{2\mu} + \frac{M^2\left(r_k+r_{k+1}\right)^2}{4\mu^2}\right). 
	\end{equation}
	From Lemma~\ref{lemma:r-up}, we have 
	\begin{equation}\label{eq:rk-up}
		r_{k+1} \leq \frac{\frac{1}{2\mu} \cdot M r_k + \sigma_k}{1-\sigma_k} \cdot r_k, \text{ if } \sigma_k<1.
	\end{equation}
	Now we prove the results of Lemma~\ref{lem:base} by induction. 
	
	1) When $k=0$, from Eq.~\eqref{eq:init-cond}, we have $\sigma_0 \leq \frac{q}{1+q} < 1$ and $\frac{M r_0}{\mu} \leq \frac{q}{1+q}-\sigma_0$. Hence, we obtain
	\[ r_{1} \stackrel{\eqref{eq:rk-up}}{\leq} \frac{\frac{1}{2\mu} \cdot M r_0 + \sigma_0}{1-\sigma_0}r_0 \stackrel{\eqref{eq:init-cond}}{\leq} \frac{\left(\frac{q}{1+q}-\sigma_0\right)+\sigma_0}{1-\frac{q}{1+q}} r_0 = q r_0. \]
	Thus, Eq.~\eqref{eq:linear-con} holds for $k=0$.
	
	Now we turn to Eq.~\eqref{eq:sigma-bound}.
	The first inequality is trivial, and we prove a stronger result for the second inequality.
	For brevity, we denote $a := M r_0 / \mu$. Then we have
	\begin{equation}\label{eq:a}
		a = \frac{M r_0}{\mu} \stackrel{\eqref{eq:init-cond}}{\leq} \frac{q}{1+q} < q.
	\end{equation}
	Moreover, when $0 < t \leq \frac{q(1-q)}{8}$, we have $ t \leq \frac{q(1-q)}{8} \leq \frac{q(1-q)}{2(1+q)^2}<1$ since $0<q<1$. Thus, we derive
	\begin{equation}\label{p1}
		\frac{1+q}{1-q} \cdot \left(t+\frac{q}{1+q}t^2\right) + \frac{t q}{(1+q)^2} \leq \frac{1+q}{1-q} \cdot \left(1+\frac{q}{1+q}\right)t + \frac{t}{1-q} = \frac{2+2q}{1-q}t \leq \frac{q}{1+q}.
	\end{equation}
	Set $b := t\left(\frac{q}{1+q}-\sigma_0\right)$. By $t<1$, we have
	\begin{equation}\label{eq:b}
	    b \leq \frac{q}{1+q}t < q.
	\end{equation}
	Next, multiplying $\left(\frac{q}{1+q}-\sigma_0\right)$ in both sides of Eq.~\eqref{p1} and noting that $b \leq \frac{q}{1+q}t$, we obtain 
	\[ \frac{1+q}{1-q} \cdot \left(b+b^2\right) + \frac{b q}{(1+q)^2} \leq \frac{q}{1+q}\left(\frac{q}{1+q}-\sigma_0\right) \stackrel{\eqref{eq:init-cond}}{\leq} \left(\frac{q}{1+q}\right)^2-\sigma_0^2, \]
	which leads to
	\[ 1 \stackrel{\eqref{eq:b}}{>} \left(\frac{q-b/2}{1+q}\right)^2 \geq \sigma_0^2 + \frac{1+q}{1-q} \cdot \left(b+b^2\right).  \]
	Hence, when $a\leq \frac{q(1-q)}{8}\left(\frac{q}{1+q}-\sigma_0\right)$, we obtain that 
	\begin{equation}\label{eq:only}
		\begin{aligned}
			1>\left(\frac{q-a/2}{1+q}\right)^2 \geq \sigma_0^2 + \frac{1+q}{1-q} \cdot \left(a+a^2\right) 
			\stackrel{\eqref{eq:a}}{\Leftrightarrow}
			0 < \frac{\frac{a}{2}+\sqrt{\sigma_0^2 + \frac{1+q}{1-q} \cdot \left(a+a^2\right)}}{1-\sqrt{\sigma_0^2 + \frac{1+q}{1-q} \cdot \left(a+a^2\right)}} \leq q.
		\end{aligned}
	\end{equation}
	Thus the second inequality of Eq.~\eqref{eq:sigma-bound} follows the first part of Eq.~\eqref{eq:only}:
	\[ \sigma_0^2 + \frac{1+q}{1-q} \cdot \left(a+a^2\right) \leq \left(\frac{q-a/2}{1+q}\right)^2 \stackrel{\eqref{eq:a}}{\leq} \left(\frac{q}{1+q}\right)^2. \] 

	2) Now let $k \geq 0$, and suppose that Eqs.~\eqref{eq:sigma-bound} and \eqref{eq:linear-con} have already been proved for all indices up to $ k$.
	Then from Eq.~\eqref{eq:sk-up} and $\sigma_i < \frac{q}{q+1} < 1, r_{i+1} \leq q r_i, \forall i \leq k$ by the inductive assumption, we could obtain 
	\begin{eqnarray}
		\sigma_{k+1}^2 & \stackrel{\eqref{eq:sk-up}}{\leq} & \sigma_k^2 + 2 \cdot \left(\frac{M \left(r_k+r_{k+1}\right)}{2 \mu} + \frac{M^2\left(r_k+r_{k+1}\right)^2}{4\mu^2}\right) \leq \dots \nonumber \\ 
		&\stackrel{\eqref{eq:linear-con}}{\leq} & \sigma_0^2 + \frac{M(1+q)}{\mu} \sum_{i=0}^{k} r_i + \frac{M^2(1+q)^2}{\mu^2} \sum_{i=0}^{k} r_i^2 \nonumber \\
		&\stackrel{\eqref{eq:linear-con}}{\leq} & \sigma_0^2 + \frac{1+q}{1-q}\left( \frac{M r_0}{\mu} + \frac{M^2 r_0^2}{\mu^2}\right) \stackrel{\eqref{eq:only}}{\leq} \left(\frac{q-a/2}{1+q}\right)^2 \stackrel{\eqref{eq:a}}{\leq} \left(\frac{q}{1+q}\right)^2<1, \label{eq:sigma-k+1}
	\end{eqnarray}
	showing that Eq.~\eqref{eq:sigma-bound} holds for $k+1$.
	Moreover, we have that 
	\begin{align*}
	    r_{k+2} &\stackrel{\eqref{eq:sigma-k+1}\eqref{eq:rk-up}}{\leq} \frac{\frac{1}{2\mu} \cdot M r_{k+1} + \sigma_{k+1}}{1-\sigma_{k+1}} \cdot r_{k+1} \stackrel{\eqref{eq:linear-con}\eqref{eq:sigma-k+1}}{\leq} \frac{\frac{a}{2}+\sqrt{\sigma_0^2 + \frac{1+q}{1-q} \cdot \left(a+a^2\right)}}{1-\sqrt{\sigma_0^2 + \frac{1+q}{1-q} \cdot \left(a+a^2\right)}} \cdot r_{k+1} 
	    \stackrel{\eqref{eq:only}}{\leq} q r_{k+1}.
	\end{align*}
    Hence, Eq.~\eqref{eq:linear-con} holds for $k+1$. We finish the inductive step. \hfill$\square$
\end{proof}

\subsection{Proof of Lemma~\ref{lem:super}}

\begin{proof}
	The results of Lemma~\ref{lem:base} still hold because the assumptions are the same.
	Using Lemma~\ref{lem:base}, we have
	\begin{equation*}
		\norm{\mI_n-\mJ_*^{-1}\mB_k} = \norm{\mJ_*^{-1}\left(\mB_k-\mJ_{*}\right)} \stackrel{\eqref{eq:jac-measure}}{\leq} \sigma_k \stackrel{\eqref{eq:sigma-bound}}{\leq} \frac{q}{1+q} < 1.
	\end{equation*}
	By Lemma~\ref{aux-lemma2}, we obtain $\mJ_*^{-1}\mB_k$ is non-singular and
	\begin{equation}\label{eq:jb}
		\norm{(\mJ_*^{-1}\mB_k)^{-1}} \stackrel{\eqref{eq:inv-i-e}}{\leq} \frac{1}{1-\norm{\mI_n-\mJ_*^{-1}\mB_k}} \leq \frac{1}{1-\sigma_k} \leq 1+q.
	\end{equation}
	Hence we obtain
	\begin{equation}\label{eq:lower}
		\frac{\norm{\mJ_*^{-1}\mB_k\vu_k}}{\norm{\vu_k}} \geq \frac{1}{\norm{(\mJ_*^{-1}\mB_k)^{-1}}} \stackrel{\eqref{eq:jb}}{\geq} \frac{1}{1+q}.
	\end{equation}
	By Lemma~\ref{lemma:sigma-up}, we have
	\begin{eqnarray*}
		\sigma_{k+1}^2 &\stackrel{\eqref{eq:sigma-update2}}{\leq} & \sigma_k^2 - \frac{f_{k+1}^2}{f_k^2} \cdot \frac{\norm{\mJ_*^{-1}\mB_k\vu_k}^2}{\norm{\vu_k}^2} + 2 \left(\sigma_k\norm{\overline{\mJ_k}} + \norm{\overline{\mJ_k}}^2\right) \\ 
		&\stackrel{\eqref{eq:lower}}{\leq}& \sigma_k^2 - \frac{f_{k+1}^2}{f_k^2} \cdot \frac{1}{(1+q)^2} + 2 \left(\sigma_k\norm{\overline{\mJ_k}} + \norm{\overline{\mJ_k}}^2\right) \\
		&\stackrel{\eqref{eq:sigma-bound}}{\leq}& \sigma_k^2 - \frac{f_{k+1}^2}{f_k^2} \cdot \frac{1}{(1+q)^2}+ 2\left(\norm{\overline{\mJ_k}} + \norm{\overline{\mJ_k}}^2\right) \\
		&\stackrel{\eqref{eq:over-j}}{\leq}& \sigma_k^2 - \frac{f_{k+1}^2}{f_k^2} \cdot \frac{1}{(1+q)^2} + \frac{M\left(r_{k+1}+r_k\right)}{\mu} + \frac{M^2\left(r_{k+1}+r_k\right)^2}{\mu^2} \\
		&\stackrel{\eqref{eq:linear-con}}{\leq}& \sigma_k^2 - \frac{f_{k+1}^2}{f_k^2} \cdot \frac{1}{(1+q)^2} +  \frac{M(1+q)r_k}{\mu} + \frac{M^2(1+q)^2r_k^2}{\mu^2}.
	\end{eqnarray*}
	Rearranging and summing up from $0$ to $k-1$,  we obtain
	\begin{eqnarray}
		\frac{1}{(1+q)^2}\sum_{i=0}^{k-1} \frac{f_{i+1}^2}{f_i^2} & \leq & \sigma_0^2 - \sigma_{k}^2 + \sum_{i=0}^{k-1} \left[ \frac{M(1+q)r_i}{\mu} + \frac{M^2(1+q)^2r_i^2}{\mu^2} \right] \nonumber \\
		&\stackrel{\eqref{eq:linear-con}}{\leq} & \sigma_0^2 + \frac{1+q}{1-q}\left( \frac{M r_0}{\mu} + \frac{M^2r_0^2}{\mu^2} \right) \stackrel{\eqref{eq:sigma-bound}}{\leq} \left(\frac{q}{1+q}\right)^2.\label{eq:theta-sum}
	\end{eqnarray}
	Hence, by the arithmetic-geometric mean inequality, we obtain
	\begin{equation*}
		\left[\frac{f_{k}}{f_0}\right]^2 = \prod_{i=0}^{k-1} \left[\frac{f_{i+1}}{f_i}\right]^2
		\leq \left[\frac{1}{k}\sum_{i=0}^{k-1} \frac{f_{i+1}^2}{f_i^2} \right]^k \stackrel{\eqref{eq:theta-sum}}{\leq} \left[\frac{q^2}{k}\right]^{k},
	\end{equation*}
	which concludes the proof.  \hfill$\square$
\end{proof}

\section{Missing Proofs in Section~\ref{subsec:bad}}

Employing $\tau_{*}(\mH)$ with other measures defined in Section~\ref{sec:preliminaries}, we can build up the relation between $\tau_{k+1}$ and $\tau_k$ as below.

\begin{lemma}\label{lemma:sigma-up-bad}
	Following Algorithm~\ref{algo:broyden-bad}, it holds that for all $k\geq 0$,
	\begin{equation}\label{eq:sigma-update1-bad}
		\tau_{k+1}^2 \leq \tau_k^2 -  \frac{\norm{ \mJ_*\left(\mJ_k^{-1} - \mH_k\right)\vy_k}^2}{\norm{\vy_k}^2} + 2 \left(\tau_k \norm{\widetilde{\mJ_k^{-1}}} + \norm{\widetilde{\mJ_k^{-1}}}^2 \right),
	\end{equation}
	if $\mJ_k$ is nondegenerate, where $\widetilde{\mJ_k^{-1}} := \mJ_*\left(\mJ_k^{-1} - \mJ_*^{-1}\right)$ follows from Eq.~\eqref{jac-measure-bad}. 
	Moreover, we can introduce $F_k$ into the above inequality:
	\begin{equation}\label{eq:sigma-update2-bad}
		\begin{aligned}
			\tau_{k+1}^2 &\leq \tau_k^2 -
			\frac{F_{k+1}^2/F_k^2}{\left(1+\norm{\left(\mJ_k-\mJ_*\right)\mJ_*^{-1}}\right)^2 \cdot \norm{\mJ_*\mH_k}^2 \cdot \norm{(\mJ_*\mH_k)^{-1}}^2} + 2 \left(\tau_k \norm{\widetilde{\mJ_k^{-1}}} + \norm{\widetilde{\mJ_k^{-1}}}^2 \right).
		\end{aligned}
	\end{equation}
\end{lemma}

\begin{proof}
	Following the definition of $\tau_{k+1}$ and $\vy_k = \mJ_k\vu_k$, we have 
	\begin{equation}\label{eq:bde-bad}
		\mJ_*\mH_{k+1} - \mI_n \stackrel{\eqref{eq:h-up}}{=} \mJ_*\mH_{k} - \mI_n +  \frac{\mJ_*\left(\mJ_k^{-1} - \mH_k\right)\vy_k\vy_k^\top}{\vy_k^\top\vy_k}.
	\end{equation}
	Taking the Frobenius norm for both sides of the above equation, we have
	\begin{eqnarray}\label{eq:part1-bad}
		&&\tau_{k+1}^2 \stackrel{\eqref{jac-measure-bad}}{=} \norm{\mJ_*\mH_{k+1} - \mI_n}_F^2 \stackrel{\eqref{eq:bde-bad}}{=} \norm{ \mJ_*\left(\mH_{k} - \mJ_*^{-1}\right) + \frac{\mJ_*\left(\mJ_k^{-1} - \mH_k\right)\vy_k\vy_k^\top}{\vy_k^\top\vy_k}}_F^2 \\ \nonumber
		&\stackrel{\eqref{jac-measure-bad}}{=}& \tau_k^2 + \frac{\norm{\mJ_*\left(\mJ_k^{-1} - \mH_k\right)\vy_k}^2}{\norm{\vy_k}^2} + 2 \cdot \frac{\vy_k^\top\left(\mH_{k}-\mJ_*^{-1}\right)^\top\mJ_*^\top\mJ_*\left(\mJ_k^{-1} - \mH_k\right)\vy_k}{\vy_k^\top\vy_k} \\ \nonumber
		&=& \tau_k^2 + \frac{\norm{\mJ_*\left(\mJ_k^{-1} - \mH_k\right)\vy_k}^2}{\norm{\vy_k}^2} + 2 \cdot \frac{\vy_k^\top\left(\mH_{k}-\mJ_k^{-1}\right)^\top \mJ_*^\top\mJ_* \left(\mJ_k^{-1} - \mH_k\right)\vy_k}{\vy_k^\top\vy_k} \\ \nonumber
		&&+ 2 \cdot \frac{\vy_k^\top\left(\mJ_k^{-1}-\mJ_*^{-1}\right)^\top \mJ_*^\top\mJ_* \left(\mJ_k^{-1} - \mH_k\right)\vy_k}{\vy_k^\top\vy_k} \\ \nonumber
		&=& \tau_k^2-\frac{\norm{\mJ_*\left(\mJ_k^{-1} - \mH_k\right)\vy_k}^2}{\norm{\vy_k}^2} + 2 \cdot \frac{\vy_k^\top\left(\mJ_k^{-1} - \mJ_*^{-1}\right)^\top\mJ_*^\top\mJ_*\left(\mJ_k^{-1} - \mH_k\right)\vy_k}{\vy_k^\top\vy_k}.
	\end{eqnarray}
	Note that the final term in Eq.~\eqref{eq:part1-bad} is
	\begin{eqnarray} \nonumber
		&&\frac{\vy_k^\top\left(\mJ_k^{-1}-\mJ_*^{-1}\right)^\top\mJ_*^\top\mJ_*\left(\mJ_k^{-1} - \mH_k\right)\vy_k}{\vy_k^\top\vy_k} \leq \frac{\norm{ \mJ_* \left(\mJ_k^{-1}-\mJ_*^{-1}\right) \vy_k}}{\|\vy_k\|}\cdot \frac{\norm{\mJ_*\left(\mJ_k^{-1} - \mH_k\right)\vy_k}}{\|\vy_k\|} \\ \nonumber
		&\leq& \norm{\mJ_* \left(\mJ_k^{-1}-\mJ_*^{-1}\right)} \cdot \left[ \norm{\mJ_*\left(\mJ_*^{-1} - \mH_k\right)} + \norm{\mJ_*\left(\mJ_k^{-1} - \mJ_*^{-1}\right)} \right] \stackrel{\eqref{jac-measure-bad}}{\leq} \left(\tau_k \norm{\widetilde{\mJ_k^{-1}}} + \norm{\widetilde{\mJ_k^{-1}}}^2\right). \label{eq:part2-bad}
	\end{eqnarray} 
	Combining Eq.~\eqref{eq:part1-bad} and Eq.~\eqref{eq:part2-bad}, we obtain Eq.~\eqref{eq:sigma-update1-bad}.

	Finally, noticing that $\vy_k=\mJ_k\vu_k$, we have
	\begin{eqnarray}\label{eq:nu-cheng} \nonumber
		\frac{F_{k+1}}{F_k} &\stackrel{\eqref{eq:fk}\eqref{eq:fjb}}{=}& \frac{\norm{\left(\mB_k-\mJ_k\right)\vu_k}}{\norm{\mB_k\vu_k}}= \frac{\norm{\left(\mH_k^{-1}-\mJ_k\right)\mJ_k^{-1}\vy_k}}{\norm{\mH_k^{-1}\mJ_k^{-1}\vy_k}} \\
		&=& \frac{\norm{\mH_k^{-1} \left(\mJ_k^{-1}-\mH_k\right)\vy_k}}{\norm{\mH_k^{-1}\mJ_k^{-1}\vy_k}} \leq \frac{\norm{\mH_k^{-1}\mJ_*^{-1}} \cdot \norm{\mJ_* \left(\mJ_k^{-1}-\mH_k\right)\vy_k}}{\norm{\mH_k^{-1}\mJ_k^{-1}\vy_k}} \nonumber \\ &\leq& \frac{\norm{(\mJ_*\mH_k)^{-1}} \cdot \norm{\mJ_* \left(\mJ_k^{-1}-\mH_k\right)\vy_k}}{\norm{\vy_k}/\norm{\mJ_k\mH_k}}
	\end{eqnarray}
	Moreover, we note that
	\begin{equation}\label{eq:jkhk}
		\norm{\mJ_k\mH_k} \leq \norm{\mJ_*\mH_k} + \norm{\left(\mJ_k-\mJ_*\right)\mH_k} \leq \norm{\mJ_*\mH_k} + \norm{\left(\mJ_k-\mJ_*\right)\mJ_*^{-1}} \cdot \norm{\mJ_*\mH_k}.
	\end{equation}
	Hence, we obtain
	\begin{equation*}
		\frac{F_{k+1}}{F_k} \stackrel{\substack{\eqref{eq:nu-cheng} \\ \eqref{eq:jkhk}}}{\leq} \left(1+\norm{\left(\mJ_k-\mJ_*\right)\mJ_*^{-1}}\right) \cdot \norm{\mJ_*\mH_k} \cdot \frac{\norm{(\mJ_*\mH_k)^{-1}} \cdot \norm{\mJ_* \left(\mJ_k^{-1}-\mH_k\right)\vy_k}}{\norm{\vy_k}}.
	\end{equation*}
	Then we introduce $ F_{k+1} / F_k$ into Eq.~\eqref{eq:sigma-update1-bad} and obtain Eq.~\eqref{eq:sigma-update2-bad}. \hfill$\square$
\end{proof}

We can also obtain the relation between $R_{k+1}$ and $R_k$ defined in Eq.~\eqref{eq:R}.
\begin{lemma}\label{lemma:r-up-bad}
	Following Algorithm~\ref{algo:broyden-bad}, we have that for all $k \geq 0$,
	\begin{equation*}
		R_{k+1} \leq \left[\tau_k + \frac{\left(1+\tau_k\right)MR_k}{2\mu^2}\right] R_k, \ \forall k\geq 0.
	\end{equation*}	
\end{lemma}
\begin{proof}
	By the update of $\vx_{k}$ in Algorithm~\ref{algo:broyden-bad}, we can obtain that
	\begin{eqnarray*}
		\mJ_*\left(\vx_{k+1}-\vx_*\right) &\stackrel{\eqref{eq:h-up}}{=}& \mJ_*\left[\vx_{k}-\mH_k\mF(\vx_k)-\vx_*\right] =  \mJ_*\left(\vx_k-\vx_*\right)-\mJ_*\mH_k\mF(\vx_k) \\
		&=& -\left[\mJ_*\mH_k-\mI_n\right]\mJ_* \left(\vx_k-\vx_*\right)-\mJ_*\mH_k \left[\mF(\vx_k)-\mJ_*\left(\vx_k-\vx_*\right)\right].
	\end{eqnarray*}	
	Taking the spectral norm in the both sides of the above equation, we have that
	\begin{eqnarray*}
		R_{k+1} &\leq& \norm{\mJ_*\mH_k-\mI_n} \cdot R_k + \norm{\mJ_*\mH_k} \cdot \norm{\mF(\vx_k)-\mJ_*\left(\vx_k-\vx_*\right)} \\ 
		&\stackrel{\eqref{jac-measure-bad}}{\leq}& \tau_k R_k + \left(1+\tau_k\right) \cdot \norm{\mF(\vx_k)-\mJ_*\left(\vx_k-\vx_*\right)} \stackrel{\eqref{eq:lemma-j2}}{\leq} \tau_k R_k  + \left(1+\tau_k\right) \cdot \frac{Mr_k^2}{2} \\ &\stackrel{\eqref{eq:rR}}{\leq}& \tau_k R_k  + \left(1+\tau_k\right) \cdot \frac{MR_k^2}{2\mu^2}.
	\end{eqnarray*}
	Therefore, we reach that
	\begin{equation*}
		R_{k+1} \leq \left[\tau_k + \frac{\left(1+\tau_k\right)MR_k}{2\mu^2}\right] R_k, \ \forall k\geq 0.
	\end{equation*} \hfill$\square$
\end{proof}

\subsection{Proof of Lemma~\ref{lem:base-bad}}
\begin{proof}
	Noting that $\mJ_k = \int_{0}^1 \mJ(\vx_k+t\vu_k)dt$, we have 
	\begin{equation*}
		\norm{\mI_n-\mJ_k\mJ_*^{-1}} = \norm{\left(\mJ_k-\mJ_*\right)\mJ_*^{-1}} \leq \norm{\mJ_k-\mJ_*} \cdot \norm{\mJ_*^{-1}} \stackrel{\eqref{eq:lemma-j1}}{\leq} \frac{M\left(r_k+r_{k+1}\right)}{2\mu} \stackrel{\eqref{eq:rR}}{\leq} \frac{M\left(R_k+R_{k+1}\right)}{2\mu^2}.
	\end{equation*}
	If $ M\left(R_k+R_{k+1}\right) \leq \mu^2$, then by Lemma~\ref{aux-lemma2}, $\mJ_k\mJ_*^{-1}$, as well as $\mJ_k$ is non-singular, and
	\begin{eqnarray}\label{eq:inv-j}
		\norm{\widetilde{\mJ_k^{-1}}} &\stackrel{\eqref{jac-measure-bad}}{=}& \norm{\left(\mJ_k\mJ_*^{-1}\right)^{-1}-\mI_n} \stackrel{\eqref{eq:inv-i-e-i}}{\leq} \frac{\norm{\mI_n-\mJ_k\mJ_*^{-1}}}{1-\norm{\mI_n-\mJ_k\mJ_*^{-1}}} \\ \nonumber
		&\leq& \frac{M\left(R_k+R_{k+1}\right)}{2\mu^2 - M\left(R_k+R_{k+1}\right)} \leq \frac{M\left(R_k+R_{k+1}\right)}{\mu^2}.
	\end{eqnarray}
	Combining with Lemma~\ref{lemma:sigma-up-bad}, we conclude that if $M\left(R_k+R_{k+1}\right) \leq \mu^2$, then 
	\begin{equation}\label{eq:sk-up-bad}
		\tau_{k+1}^2 \stackrel{\eqref{eq:sigma-update1-bad}}{\leq} \tau_k^2 + 2  \left(\tau_k \norm{\widetilde{\mJ_k^{-1}}} + \norm{\widetilde{\mJ_k^{-1}}}^2 \right) \stackrel{\eqref{eq:inv-j}}{\leq} \tau_k^2 + 2 \left[\frac{\tau_k  M\left(R_k+R_{k+1}\right)}{\mu^2} + \frac{M^2\left(R_k+R_{k+1}\right)^2}{\mu^4} \right].
	\end{equation}
	Moreover, from Lemma~\ref{lemma:r-up-bad}, we have 
	\begin{equation}\label{eq:rk-up-bad}
		R_{k+1} \leq \left[\tau_k + \frac{\left(1+\tau_k\right)MR_k}{2\mu^2}\right] R_k.
	\end{equation}
	Now	we prove by induction. 
	
	1) When $k=0$, from Eq.~\eqref{eq:init-cond-bad}, we have $\tau_0 \leq q<1$, $\frac{MR_0}{\mu^2} \leq q-\tau_0$. Then we have
	\begin{eqnarray*}
		R_{1} &\stackrel{\eqref{eq:rk-up-bad}}{\leq}& \left[\tau_0 + \frac{\left(1+\tau_0\right)MR_0}{2\mu^2}\right] R_0 \leq \left(\tau_0 + \frac{M R_0}{\mu^2}\right) R_0 \stackrel{\eqref{eq:init-cond-bad}}{\leq} \left[\tau_0 + \left(q-\tau_0\right)\right]R_0 = qR_0.
	\end{eqnarray*}
	Thus, Eq.~\eqref{eq:linear-con-bad} holds for $k = 0$.
    
    Now we turn to Eq.~\eqref{eq:sigma-bound-bad}. The first inequality is trivial, and we prove a stronger result for the second inequality.
    For brevity, we denote $A:= MR_0/\mu^2$. Then we have
	\begin{equation}\label{eq:a-bad}
		A \stackrel{\eqref{eq:init-cond-bad}}{\leq} q.
	\end{equation}
	When $t\leq \frac{q(1-q)}{6} \leq \frac{1}{4}$, we have
	\begin{equation}\label{p1-bad}
		\frac{1+q}{1-q} \cdot \left(2t+q \cdot 2t^2\right) + 2qt \leq \frac{1+q}{1-q} \cdot \left(2+\frac{q}{2}\right)t + \frac{t}{2(1-q)} \leq \frac{6t}{1-q} \leq q.
	\end{equation}
	Setting $B := t\left(q-\tau_0\right) \leq qt$, and multiplying $\left(q-\tau_0\right)$ in both sides of Eq.~\eqref{p1-bad}, we obtain 
	\[ \frac{1+q}{1-q} \cdot \left(2B+2B^2\right) + 2q B\leq q(q-\tau_0) \stackrel{\eqref{eq:init-cond-bad}}{\leq} q^2-\tau_0^2,
	\]
	which leads to
	\[ \left(q-B\right)^2 \geq \tau_0^2 + \frac{1+q}{1-q} \cdot \left(2B+2B^2\right).
	\]
	Hence when $A\leq \frac{q(1-q)}{6}\cdot \left(q-\tau_0\right)$, we obtain 
	\begin{equation}\label{eq:only-bad}
		\begin{aligned}
			\left(q-A\right)^2 \geq \tau_0^2 + \frac{1+q}{1-q} \cdot \left(2A+2A^2\right)
			\stackrel{\eqref{eq:a-bad}}{\Leftrightarrow}
			\sqrt{\tau_0^2 + \frac{1+q}{1-q} \cdot \left(2A+2A^2\right)} + A \leq q,
		\end{aligned}
	\end{equation}
	Thus the second inequality of Eq.~\eqref{eq:sigma-bound-bad} follows the first part of Eq.~\eqref{eq:only-bad}:
	\[ \tau_0^2 + \frac{1+q}{1-q} \cdot \left(2A+2A^2\right) \leq \left(q-A\right)^2 \stackrel{\eqref{eq:a-bad}}{\leq} q^2. \]
	2) Now let $k \geq 0$, and suppose that Eqs.~\eqref{eq:sigma-bound-bad} and \eqref{eq:linear-con-bad} hold for all indices up to $ k$. 
	Then it holds that $\frac{M (R_k+R_{k+1})}{2\mu^2} \stackrel{\eqref{eq:linear-con-bad}}{\leq} \frac{M R_0}{\mu^2} \stackrel{\eqref{eq:init-cond-bad}}{\leq} \frac{1}{2}.
	$
	Therefore, Eq.~\eqref{eq:sk-up-bad} holds for the index $k$. 
	Noticing that $\tau_i \stackrel{\eqref{eq:sigma-bound-bad}}{\leq } q \leq 1, R_{i+1} \leq q R_i, \forall i \leq k$, then by the inductive assumption, we could obtain
	\begin{eqnarray}\nonumber
		\tau_{k+1}^2 &\stackrel{\eqref{eq:sk-up-bad}}{\leq}& \tau_k^2 + \frac{2M\left(R_k+R_{k+1}\right)}{\mu^2}+\frac{2M^2\left(R_k+R_{k+1}\right)^2}{\mu^4} \leq ... \\ \nonumber
		&\stackrel{\eqref{eq:linear-con-bad}}{\leq}& \tau_0^2 + \frac{2M(1+q)}{\mu^2} \cdot \sum_{i=0}^{k} R_i + \frac{2 M^2(1+q)^2}{\mu^4} \cdot \sum_{i=0}^{k} R_i^2 \\ 
		&\stackrel{\eqref{eq:linear-con-bad}}{\leq}& \tau_0^2 + \frac{1+q}{1-q}\left(\frac{2 M R_0}{\mu^2}+\frac{2 M^2 R_0^2}{\mu^4}\right) \stackrel{\eqref{eq:only-bad}}{\leq} \left(q-A\right)^2 \stackrel{\eqref{eq:a-bad}}{\leq} q^2, \label{eq:tau-k+1}
	\end{eqnarray}
	showing that Eq.~\eqref{eq:sigma-bound-bad} holds for $k+1$.
	Moreover, we have that
	\begin{eqnarray*}
	    R_{k+2} &\stackrel{\eqref{eq:rk-up-bad}}{\leq}& \left[\tau_{k+1} + \frac{\left(1+\tau_{k+1}\right)MR_{k+1}}{2\mu^2}\right]R_{k+1} \\ &\stackrel{\eqref{eq:linear-con-bad}\eqref{eq:tau-k+1}}{\leq}& \left[ \sqrt{\tau_0^2 + \frac{1+q}{1-q} \cdot \left(2A+2A^2\right)}+A\right]R_{k+1} \stackrel{\eqref{eq:only-bad}}{\leq} q R_{k+1}.
	\end{eqnarray*}
	Hence, Eq.~\eqref{eq:linear-con-bad} holds for $k+1$. We finish the inductive step. \hfill$\square$
\end{proof}

\subsection{Proof of Lemma~\ref{lem:super-bad}}
\begin{proof}
	The results of Lemma~\ref{lem:base-bad} still hold because the assumptions are the same.
	Using Lemma~\ref{lem:base-bad} and $q \leq 1/2$, we have
	\begin{eqnarray}\label{eq:s1}
		\norm{\left(\mJ_k-\mJ_*\right)\mJ_*^{-1}} &\leq& \norm{\left(\mJ_k-\mJ_*\right)} \cdot \norm{\mJ_*^{-1}} \stackrel{\eqref{eq:lemma-j1}}{\leq} \frac{M (r_k+r_{k+1})}{2\mu} \stackrel{\eqref{eq:rR}}{\leq} \frac{M (R_k+R_{k+1})}{2\mu^2} \nonumber \\
		&\stackrel{\eqref{eq:linear-con-bad}}{\leq}& \frac{M R_0}{\mu^2} \stackrel{\eqref{eq:init-cond-bad}}{\leq} \frac{q(1-q)}{6} \leq \frac{1}{24},
	\end{eqnarray}
	and 
	\begin{equation*}
		\norm{\mI_n-\mJ_*\mH_k}_F \stackrel{\eqref{jac-measure-bad}}{=} \tau_k \stackrel{\eqref{eq:sigma-bound-bad}}{\leq} q \leq 1/2.
	\end{equation*}
	Then
	\begin{equation}\label{eq:s2}
		\norm{\mJ_*\mH_k} \leq \norm{\mI_n} + \norm{\mI_n-\mJ_*\mH_k} \leq 1+q\leq 1.5,
	\end{equation}
	and from Lemma~\ref{aux-lemma2}, we obtain 
	\begin{equation}\label{eq:s3}
		\norm{(\mJ_*\mH_k)^{-1}} \stackrel{\eqref{eq:inv-i-e}}{\leq} \frac{1}{1-\norm{\mI_n-\mJ_*\mH_k}} \leq 2.
	\end{equation}
	Hence by Eqs.~\eqref{eq:s1},\eqref{eq:s2},\eqref{eq:s3}, we obtain 
	\begin{equation}\label{eq:multi}
	    \left(1+\norm{\left(\mJ_k-\mJ_*\right)\mJ_*^{-1}}\right)^2 \cdot \norm{\mJ_*\mH_k}^2 \cdot \norm{(\mJ_*\mH_k)^{-1}}^2 
	   \leq \left(\frac{25}{24} \times 1.5 \times 2\right)^2 \leq 10.
	\end{equation}
	From Lemma~\ref{lemma:sigma-up-bad} and Eq.~\eqref{eq:multi}, we have
	\begin{eqnarray*}
		\tau_{k+1}^2 & \stackrel{\eqref{eq:sigma-update2-bad}}{\leq}& \tau_k^2 - \frac{1}{10} \cdot \frac{F_{k+1}^2}{F_k^2} + 2  \left(\tau_k \norm{\widetilde{\mJ_k^{-1}}}+\norm{\widetilde{\mJ_k^{-1}}}^2 \right) \\ \nonumber
		&\stackrel{\eqref{eq:sigma-bound-bad}\eqref{eq:inv-j}}{\leq}& \tau_k^2 - \frac{1}{10} \cdot \frac{F_{k+1}^2}{F_k^2} + \frac{2 M\left(R_{k+1}+R_k\right)}{\mu^2} + \frac{2M^2\left(R_{k+1}+R_k\right)^2}{\mu^4} \\ \nonumber
		&\stackrel{\eqref{eq:linear-con-bad}}{\leq}& \tau_k^2 - \frac{1}{10} \cdot \frac{F_{k+1}^2}{F_k^2} + \frac{2 M(1+q)R_k}{\mu^2} + \frac{2M^2(1+q)^2R_k^2}{\mu^4}.
	\end{eqnarray*}
	Rearranging and summing up from $0$ to $k-1$, we get
	\begin{eqnarray}
		\frac{1}{10}\sum_{i=0}^{k-1} \frac{F_{i+1}^2}{F_i^2} &\leq& \tau_0^2 - \tau_{k}^2 + \sum_{i=0}^{k-1} \left[\frac{2 M (1+q)R_i}{\mu^2} + \frac{2 M^2(1+q)^2R_i^2}{\mu^4}\right] \nonumber \\ \nonumber
		& \stackrel{\eqref{eq:linear-con-bad}}{\leq}& \tau_0^2 + \frac{1+q}{1-q}\left(\frac{2 M R_0}{\mu^2} + \frac{2M^2R_0^2}{\mu^4}\right) \stackrel{\eqref{eq:sigma-bound-bad}}{\leq} q^2. \label{eq:theta-sum-bad}
	\end{eqnarray}
	Hence, by the arithmetic-geometric mean inequality, we obtain the final result
	\begin{equation}\label{eq:f-final}
		\left[\frac{F_{k}}{F_0}\right]^2 = \prod_{i=0}^{k-1} \left[\frac{F_{i+1}}{F_i}\right]^2
		\leq \left[\frac{1}{k}\sum_{i=0}^{k-1} \frac{F_{i+1}^2}{F_i^2} \right]^k \stackrel{\eqref{eq:theta-sum-bad}}{\leq} \left[\frac{10 q^2}{k}\right]^{k}.
	\end{equation}
	Finally, noting that for all $k \geq 1$,
	\begin{equation}\label{eq:lam-f}
	    \frac{F_k}{L} = \frac{F_k}{\norm{\mJ_*}} \leq f_k \stackrel{\eqref{eq:lam}}{=} \norm{ \mJ_*^{-1}\mF(\vx_k)} \leq \norm{\mJ_*^{-1}} \cdot F_k = \frac{F_k}{\mu},
	\end{equation}
	we obtain 
	\[ f_{k} \stackrel{\eqref{eq:lam-f}}{\leq} \frac{F_k}{\mu} \stackrel{\eqref{eq:f-final}}{\leq} \frac{1}{\mu} \cdot \left[\frac{10 q^2}{k}\right]^{k/2} F_0 \stackrel{\eqref{eq:lam-f}}{\leq} \frac{L}{\mu} \cdot \left[\frac{10 q^2}{k}\right]^{k/2} f_0 = \varkappa \left[\frac{10 q^2}{k}\right]^{k/2} f_0. \]
	The proof is finished. \hfill$\square$
\end{proof}

\subsection{Proof of Theorem~\ref{thm:bad}}\label{app:thmbad}
\begin{proof}
	To give an upper and lower bound of $q_m$, we consider 
	\[ f(q) := q(1-q)\left(q-\tau_0\right), \ q \geq \tau_0. \]
	We only need to consider the second constraint in the definition of $q_m$.
	
	Note that $q(1-q)$ is increasing when $q \leq \frac{1}{2}$, and $q$ is also 
	increasing when $q>0$. Hence $f(q)$ is increasing when $\frac{1}{2} \geq q \geq \tau_0$, where we use $\frac{1}{2} \geq  \tau_0$ by Eq.~\eqref{eq:init-super-bad}.  Moreover, we have
    \[ f\left(\frac{1}{2}\right) = \frac{1}{4} \left(\frac{1}{2}-\tau_0\right) \stackrel{\eqref{eq:init-super-bad}}{\geq} \frac{6M R_0}{\mu^2}, \]
	showing that $q=\frac{1}{2}$ is in the feasible set. Thus $q_m$ is well-defined and
	\begin{equation}\label{eq:qm-2-bad}
	    0 \leq \tau_0 \leq q_m \leq \frac{1}{2}.
	\end{equation}
	Since $f(\tau_0) = 0 \leq \frac{6M R_0}{\mu^2} \leq f(\frac{1}{2})$, we get $f(q_m) = \frac{6M R_0}{\mu^2}$. 
	Then we obtain
    \[ \frac{6 M R_0}{\mu^2} = f(q_m) \stackrel{\eqref{eq:qm-2-bad}}{\geq} \frac{q_m}{2}\left(q_m-\tau_0\right), q_m \geq \tau_0, \]
	which gives
	\begin{equation}\label{eq:qm-upper-bad}
		q_m \stackrel{\eqref{eq:qm-2-bad}}{\leq} \frac{\tau_0}{2} + \frac{1}{2} \sqrt{\tau_0^2+\frac{48 M R_0}{\mu^2}} \leq \tau_0+4\sqrt{\frac{ M R_0}{\mu^2}} \leq 4 \left(\tau_0+\sqrt{\frac{ M R_0}{\mu^2}}\right).
	\end{equation}
	Furthermore, we have that
	\[ \frac{6 M R_0}{\mu^2} = f(q_m) \leq q_m\left(q_m-\tau_0\right), q_m \geq \sigma_0, \]
	which gives
	\begin{equation}\label{eq:qm-lower-bad}
	    q_m \stackrel{\eqref{eq:qm-2-bad}}{\geq} \frac{\tau_0}{2} +\frac{1}{2}\sqrt{\tau_0^2+\frac{24 M R_0}{\mu^2}} \geq \frac{\tau_0}{2} + 2\sqrt{\frac{M R_0}{\mu^2}} \geq \frac{1}{2} \left(\tau_0+\sqrt{\frac{ M R_0}{\mu^2}}\right).
	\end{equation}
	Therefore, combining Eqs.~\eqref{eq:qm-upper-bad} and \eqref{eq:qm-lower-bad}, we obtain $q_m = \Theta\left(\tau_0+\sqrt{\frac{ M R_0}{\mu^2}}\right)$.
	Finally, Eq.~\eqref{eq:easyread-convergence-bad} can be derived by choosing $q=q_m$ in Lemmas~\ref{lem:base-bad} and~\ref{lem:super-bad}, and
	\[ f_{k} \stackrel{\eqref{eq:qm-upper-bad}}{\leq}  \left[\frac{10\times 16\left(\tau_0+\sqrt{\frac{M R_0}{\mu^2}}\right)^2}{k}\right]^{k/2}f_0 \leq \left[\frac{13\left(\tau_0+\sqrt{\frac{M R_0}{\mu^2}}\right)}{\sqrt{k}}\right]^{k}f_0. \] \hfill$\square$
\end{proof}

\bibliography{reference}
\bibliographystyle{plainnat}

\end{document}